\documentclass[11pt]{amsart}

\usepackage[T1]{fontenc}
\usepackage{microtype}

\usepackage[margin=1in]{geometry}
\usepackage{graphicx}
\graphicspath{{figures/}}
\usepackage{amsmath,amssymb,amsthm}
\usepackage{newpxtext,newpxmath}
\usepackage[numbers,sort&compress]{natbib}
\usepackage[colorlinks=true,allcolors=blue]{hyperref}
\usepackage{doi}
\usepackage{tikz}
\usetikzlibrary{arrows.meta,positioning,calc,fit,backgrounds,matrix}
\usepackage{paperA-figures}

\newcommand{\Species}{\mathcal{S}}
\newcommand{\Reactions}{\mathcal{R}}

\newcommand{\Bonds}{\mathcal{B}}
\newcommand{\Atoms}{\mathcal{T}}
\newcommand{\Centers}{\mathcal{C}}

\providecommand{\Nat}{\mathbb{N}}

\providecommand{\addbibresource}[1]{}
\newcommand{\GFspecies}{F}
\newcommand{\GFensemble}{Z}
\newcommand{\GFcomb}{\Psi}
\newcommand{\dd}{\,\mathrm{d}}
\newcommand{\Deriv}[2]{\frac{\mathrm{d}#1}{\mathrm{d}#2}}
\newcommand{\PDeriv}[2]{\frac{\partial #1}{\partial #2}}

\theoremstyle{definition}
\newtheorem{definition}{Definition}

\theoremstyle{plain}

\newtheorem{proposition}{Proposition}
\newtheorem{theorem}{Theorem}

\theoremstyle{remark}
\newtheorem*{remark}{Remark}

\title{Equilibrium Combinatorial Self-Assembly \\ via Generating Functions}
\author{Andr\'es Ortiz-Mu\~{n}oz}
\address{Harvard Medical School}
\email{andortiz19@gmail.com}
\date{}

\begin{document}

\begin{abstract}
We develop a generating-function calculus for equilibrium combinatorial self-assembly.
Starting from a bond-level specification of allowable interactions, we define a symmetry-weighted \emph{species generating function} whose evaluation yields equilibrium concentrations, and an \emph{ensemble generating function} (an exponential transform) that packages equilibrium probabilities of mixtures.
We ground these objects in both deterministic (coagulation--fragmentation) and stochastic (master equation) dynamics, showing how detailed balance leads to the equilibrium expressions and how the exponential generating function arises as a partition function.
We develop a formal power series calculus---derivatives, integrals, exponentials, composition---where each operation acquires a precise combinatorial interpretation.
The paper is organized around two regimes: \emph{cycle-free assembly}, where binding equations for the species generating function are nonlinear and the ensemble equation couples to the species generating function; and \emph{assembly with cycles}, where a cycle-opening term enters the species equation and the exponential transform linearizes the ensemble equation into a closed linear PDE with operator-exponential solutions.
Each regime is developed with a linear polymer worked example in which we compute equilibrium concentrations, extract canonical partition functions, and derive canonical factorial moments.
A cross-linked polymer example---combining heterotypic chain bonds with homotypic cross-links---illustrates both regimes together, yielding a factorized canonical partition function and an explicit gelation surface.
\end{abstract}

\maketitle


\section{Introduction}
\label{sec:intro}

Equilibrium self-assembly sits at an interface between local chemistry and global combinatorics.
A specification of admissible interactions---which sites can bind, with what affinities, and under what constraints---may be simple, yet the induced space of possible assemblies can be enormous.
At equilibrium, this space is not merely an enumeration problem: it determines which structures are abundant and how likely particular motifs (notably \emph{cycles}) are to occur.
The present approach leans on generating functions and symmetry-weighted counting in the spirit of analytic combinatorics and species \cite{flajolet2009analytic,joyal1981theorie,bergeron1998combinatorial}.

The central idea is to encode each assembly as a monomial in formal variables that track agents, bond types, and free sites, weighted by a symmetry factor.
Summing these monomials produces a \emph{species generating function} whose evaluation yields equilibrium concentrations; exponentiating it produces an \emph{ensemble generating function} encoding equilibrium distributions.
Differentiation, integration, and exponential operators then acquire precise combinatorial meanings (marking, discounting, and forming ensembles) that translate directly into equilibrium identities and recursions.
On the probabilistic side, this matches the well-known product-form structure of stochastic equilibria and connects naturally to the stochastic mechanics program \cite{whittle1986systems,baez2018quantum,baez2015quantum}.

Two themes drive the narrative.
\begin{itemize}
  \item \textbf{From structure to equilibrium.} We fix a structural language (bond systems, bond structures, validity) that cleanly separates \emph{what assemblies exist} from \emph{how they are produced}. Under detailed-balance assumptions, equilibrium weights depend only on energetic parameters and internal symmetries, making generating functions a natural summary object.
  \item \textbf{Cycle-free versus assembly with cycles.} Acyclic assemblies already exhibit rich recursion and closed-form solutions; allowing cycles introduces additional symmetry, logarithmic terms, and a genuine distinction between inter- and intramolecular binding. In the dimensionless variables, the species and ensemble generating functions retain the same $\gamma$-free form as in the acyclic case, but the ``bond discounting'' identity acquires an additional \emph{cycle-opening} term. It is precisely this term that allows the exponential change of variables to \emph{linearize} the ensemble equation---a feature specific to the cyclic regime.
\end{itemize}

\paragraph{Contributions.}
(i)~A bond-system formalism for assemblies connected to a symmetry-weighted species generating function, including the Euler characteristic as a topological admissibility criterion.
(ii)~A detailed account of both deterministic and stochastic dynamics, showing how detailed balance leads to the equilibrium expressions and the ensemble generating function arises as a partition function.
(iii)~A calculus of formal power series operations (derivatives, integrals, exponentials, composition, pointing) with combinatorial interpretations.
(iv)~For \emph{cycle-free assembly}: general binding equations, recursion by integration, and a linear polymer worked example.
(v)~For \emph{canonical observables}: composition-constrained partition functions and factorial moments from probability generating functions.
(vi)~A saddle-point approximation bridging the stochastic and deterministic descriptions, with an explicit error expansion and rigorous validity conditions (H-admissibility).
(vii)~For \emph{assembly with cycles}: the cyclic discounting identity and its linearization via exponentiation, yielding operator-exponential solutions, with a linear-polymer-with-rings worked example.
(viii)~A \emph{cross-linked polymer} worked example combining heterotypic chain bonds with homotypic cross-links, illustrating the site-derivative system, cubic implicit equation, operator-exponential factorization of canonical partition functions, and gelation analysis in both acyclic and cyclic regimes.

\paragraph{Roadmap.}
Section~\ref{sec:formalism} introduces the structural formalism.
Section~\ref{sec:dynamics} develops deterministic and stochastic dynamics and establishes the equilibrium expressions.
Section~\ref{sec:gf} collects the formal power series operations, shared canonical observables, and the saddle-point bridge between stochastic and deterministic descriptions.
Section~\ref{sec:cyclefree} treats cycle-free assembly (general theory and linear polymer example).
Section~\ref{sec:cycles} treats assembly with cycles (general theory and linear polymer example).
Section~\ref{sec:crosslinked} works through a cross-linked polymer system in both regimes, deriving gelation thresholds.
Section~\ref{sec:outlook} closes with outlook and open directions.
Appendices collect closure-operator expansions, proofs of canonical moment identities, the saddle-point derivation and H-admissibility conditions, and foundational material on species, groupoid cardinality, and stuff types.

\section{Bond systems, structures, and assemblies}
\label{sec:formalism}

This section fixes the structural objects that the rest of the paper counts, manipulates, and weights at equilibrium.

\subsection{Bond systems}

\begin{definition}
A \emph{bond system} is a tuple $(B,C,s,t,1,\cdot)$, where $B$ is a set of \emph{bonds}, $C$ is a set of \emph{content} relations, and $s,t:C\to B$ are source and target functions.
For $c\in C$ with $s(c)=x$ and $t(c)=y$, we write $x\xrightarrow{c}y$, indicating that $x$ \emph{contains} $y$.
Each bond $x\in B$ has an identity containment $1_x:x\to x$.
For $c:x\to y$ and $c':y\to z$, their composition is $cc':x\to z$.
These data must satisfy: (1)~identity laws, (2)~associativity, (3)~anti-symmetry (if $c:x\to y$ and $c':y\to x$ then $x=y$), and (4)~left cancellation.
\end{definition}

In categorical terms, a bond system is a skeletal category where all morphisms are monomorphisms.
The \emph{order} of a bond is defined hierarchically: order~0 bonds (\emph{atoms}) contain only themselves; inductively, a bond has order $n$ if the maximum order among its contained bonds is $n-1$.

\subsection{Bond structures and valid assemblies}

\begin{definition}
Given a bond system, a \emph{bond structure} $F$ assigns to each bond $b\in B$ a set $F[b]$ of instances, and to each containment $c:x\to y$ an injective function $F[c]:F[x]\to F[y]$.
\end{definition}

A bond structure is \emph{connected} if its underlying undirected graph is connected.
Physical validity is determined by which structures are reachable through assembly dynamics.
Under detailed balance, the equilibrium distribution depends only on which structures are reachable and the energy content of the bonds.

\subsection{Assembly systems}
\label{sec:assembly-systems-data}

An \emph{assembly system} consists of:
a class $\Species$ of admissible species (connected assemblies modulo symmetry),
a set $\Reactions$ of binding/unbinding events,
sets $\Atoms$ (atom types) and $\Bonds$ (bond types),
a set $\Centers$ of sites/centers (each bond consumes two compatible centers),
content maps recording how many atoms/bonds/centers each species contains,
and a stoichiometry assignment for reactions.
Each species $s\in\Species$ carries a nonneg.\ integer \emph{cycle rank} $\vartheta(s)$ counting independent intramolecular closures.

\subsection{Euler characteristic}
\label{sec:euler-characteristic}

For a bond structure $s$ with $n$ atoms (vertices) and $m$ bonds (edges), the Euler characteristic is $\chi(s) = n - m$.
Every species is by definition a \emph{connected} assembly, so its zeroth Betti number is $\beta_0(s) = 1$.
The first Betti number $\beta_1(s)$ counts independent cycles and equals the cycle rank: $\beta_1(s) = \vartheta(s)$.
The Euler characteristic identity for a connected graph then reads
\begin{equation}\label{eq:euler_identity}
\chi(s) = \beta_0(s) - \beta_1(s) = 1 - \vartheta(s),
\end{equation}
equivalently $\vartheta(s) = 1 - \chi(s) = m - n + 1$.
Trees satisfy $\chi=1$ (no cycles), unicyclic structures $\chi=0$, and $k$-cyclic structures $\chi = 1-k < 0$.
This identity is the topological foundation for the $\gamma$-power counting developed in Section~\ref{sec:deterministic}.

\section{Dynamics and equilibrium}
\label{sec:dynamics}

The structural formalism of Section~\ref{sec:formalism} specifies \emph{which} assemblies are admissible.
To make equilibrium predictions we also need an account of \emph{how} assemblies form and break.

\subsection{Deterministic dynamics: coagulation--fragmentation}
\label{sec:deterministic}

At a deterministic level, the concentrations $c_s$ of species evolve according to rate equations.
Each reaction $r\in\Reactions$ has a stoichiometry vector $\nu_r$ and a net flux $J_r(c)$.
For reversible binding $X_r + Y_r \rightleftarrows Z_r$:
\begin{equation}\label{eq:coag_frag}
\dot{c}_s = \sum_{r\in\Reactions} \nu_{r,s}\,J_r(c),
\qquad J_r(c) = \kappa_r^+\,c_{X_r}\,c_{Y_r} - \kappa_r^-\,c_{Z_r}.
\end{equation}

\paragraph{Detailed balance.}
A much stronger condition than $\dot{c}_s=0$ for all $s$ is \emph{detailed balance}: each reaction is individually balanced,
\begin{equation}\label{eq:detailed_balance}
\kappa_r^+\,c_{X_r}\,c_{Y_r} = \kappa_r^-\,c_{Z_r}
\qquad\text{for every } r\in\Reactions.
\end{equation}
This splits the infinite-dimensional system into independent equations per reaction.
Since the affinity $y_r = \kappa_r^+/\kappa_r^-$ is determined by the free-energy difference for each binding step, one can solve by induction on assembly size.
For a species $s$ built from atoms of types $a\in\Atoms$ and bonds of types $b\in\Bonds$, the equilibrium concentration is
\begin{equation}\label{eq:equilibrium_concentration}
c_s = \frac{1}{\varphi(s)}\,
\gamma^{\vartheta(s)}\,
\prod_{a\in\Atoms} c_a^{\alpha(s,a)}\,
\prod_{b\in\Bonds} y_b^{\beta(s,b)},
\end{equation}
where $\varphi(s)$ is the number of automorphisms of $s$ and $\vartheta(s)$ is the cycle rank.
In the acyclic case $\vartheta(s)=0$, so~\eqref{eq:equilibrium_concentration} reduces to the usual tree-level formula.

\paragraph{Euler characteristic and $\gamma$-power counting.}
The $\gamma$-exponent in~\eqref{eq:equilibrium_concentration} is $\vartheta(s)$, while the atom and bond factors contribute $+n$ and $-m$ powers of concentration (since $y_b$ has units of $[\text{concentration}]^{-1}$).
The total power of $\gamma$ in $c_s$ is therefore
\[
\vartheta(s) + n - m.
\]
Because $s$ is a connected species, the Euler characteristic identity~\eqref{eq:euler_identity} gives $\chi(s) = n - m = 1 - \vartheta(s)$, so
\[
\vartheta(s) + n - m = \vartheta(s) + (1 - \vartheta(s)) = 1,
\]
regardless of the topology of $s$.
Hence $c_s$ always carries exactly one power of $\gamma$, for trees and cyclic species alike.
The role of $\vartheta(s)$ is not to make individual terms consistent, but to ensure the cancellation is exact: it is the first Betti number that compensates the deficit between atom and bond counts enforced by the Euler characteristic.

\paragraph{Dimensionless normalization.}
Introduce dimensionless variables by scaling with $\gamma$:
\[
c_a = \gamma z_a,\qquad z_b = \gamma y_b
\quad\Longleftrightarrow\quad
y_b = \frac{z_b}{\gamma}.
\]
Substituting into~\eqref{eq:equilibrium_concentration} gives
\[
c_s
= \frac{1}{\varphi(s)}\,
\gamma^{\vartheta(s)}\,
\prod_{a\in\Atoms}(\gamma z_a)^{\alpha(s,a)}
\prod_{b\in\Bonds}\!\left(\frac{z_b}{\gamma}\right)^{\beta(s,b)}
= \gamma\,\frac{1}{\varphi(s)}\,
\prod_{a\in\Atoms} z_a^{\alpha(s,a)}
\prod_{b\in\Bonds} z_b^{\beta(s,b)},
\]
since the Euler characteristic identity gives $\vartheta(s) + n - m = 1$, i.e.\ $\vartheta(s) + \sum_a\alpha(s,a) - \sum_b\beta(s,b) = 1$.
Equivalently, the reduced concentration $\tilde c_s := c_s/\gamma$ is purely combinatorial:
\[
\tilde c_s
= \frac{1}{\varphi(s)}\,
\prod_{a\in\Atoms} z_a^{\alpha(s,a)}
\prod_{b\in\Bonds} z_b^{\beta(s,b)}.
\]

\subsection{Stochastic dynamics: master equation and product-form equilibrium}
\label{sec:stochastic}

A \emph{state} is a function $x:\Species\to\Nat$ with finite support.
The system evolves as a continuous-time Markov chain with the chemical master equation:
\begin{equation}\label{eq:master_equation}
\dot{p}_t(x) = \sum_{r\in\Reactions}
\Bigl[
  w_r(x - \nu_r)\,p_t(x-\nu_r) - w_r(x)\,p_t(x)
\Bigr].
\end{equation}
In the deterministic description, the bimolecular rate constant $\kappa_r^+$ has units $[\text{conc}^{-1}\!\cdot\text{time}^{-1}]$.
In the stochastic description the propensities are
\[
w_r^+(x) = \frac{\kappa_r^+}{N_A V}\,x(X_r)\,x(Y_r),
\qquad
w_r^-(x) = \kappa_r^-\,x(Z_r),
\]
so the stochastic bimolecular rate constant is $\kappa_r^+/(N_A V)$ and the unimolecular rate $\kappa_r^-$ carries no volume factor.
At stochastic detailed balance,
\[
\frac{x(Z_r)}{x(X_r)\,x(Y_r)} = \frac{\kappa_r^+}{\kappa_r^-\cdot N_A V} = \frac{y_r}{N_A V}.
\]
Converting to molar concentrations $c = x/(N_A V)$ recovers the deterministic ratio $c_{Z_r}/(c_{X_r}\,c_{Y_r}) = y_r$: volume cancels for bimolecular binding.
For an \emph{intramolecular} closure, however, both sites are on the same assembly, so the forward step is unimolecular and its propensity scales as $x$ rather than $x^2/(N_A V)$.
The equilibrium ratio then retains a factor of $1/(N_A V)$, and the concentration genuinely depends on the system volume.
This asymmetry---bimolecular binding is volume-independent, intramolecular closure is not---is the physical origin of the cycle counter $\gamma = 1/(N_A V)$.

\paragraph{Combinatorial generating function $\GFcomb$.}
The Euler characteristic identity established in Section~\ref{sec:euler-characteristic} guarantees $c_s = \gamma\,\tilde c_s$ for every species $s$, where $\tilde c_s$ is purely combinatorial in the dimensionless variables $z_a, z_b$.
Define the \emph{combinatorial species generating function}
\begin{equation}\label{eq:psi_def}
\GFcomb := \sum_{s\in\Species} \tilde c_s.
\end{equation}
Since $\lambda_s = c_s/\gamma$, the Euler characteristic identity immediately gives
\begin{equation}\label{eq:lambda_tilde_c}
\lambda_s = \frac{c_s}{\gamma} = \tilde c_s.
\end{equation}
The Poisson activity of each species equals its combinatorial weight.
Consequently $\sum_{s}\lambda_s = \GFcomb$, and $\GFcomb$ is purely combinatorial: it depends on $z_a$ and $z_b$ but carries no explicit $\gamma$.

\paragraph{Product-form equilibrium.}
Under detailed balance, the stationary distribution takes a product form:
\begin{equation}\label{eq:product_form}
p_{\mathrm{eq}}(x) = \frac{1}{\mathcal{Z}}\,\prod_{s\in\Species}
\frac{\lambda_s^{x(s)}}{x(s)!},
\end{equation}
where $\lambda_s$ is the equilibrium activity of species $s$ and $\mathcal{Z}$ is the partition function.
Using~\eqref{eq:lambda_tilde_c}, $\sum_s \lambda_s = \GFcomb$, so the partition function is
\begin{equation}\label{eq:partition_function}
\mathcal{Z} = \exp\!\left(\sum_{s\in\Species}\lambda_s\right) = \exp(\GFcomb),
\end{equation}
and the ensemble generating function is
\begin{equation}\label{eq:ensemble_gf_dynamics}
\GFensemble = \exp(\GFcomb).
\end{equation}
Both are purely combinatorial: $\GFcomb$ carries no explicit $\gamma$, because the Euler characteristic identity ensures all factors of $\gamma$ cancel exactly.
The free energy of the assembly system is $\log\mathcal{Z} = \GFcomb$, a formal power series in $z_a$ and $z_b$ that enumerates species weighted by their combinatorial symmetry factors $1/\varphi(s)$.
This is the central object of the generating-function approach developed in Section~\ref{sec:gf}.

\section{Generating functions and formal power series}
\label{sec:gf}

Let $R[\![z_1,z_2,\dots]\!]$ be the ring of formal power series in indeterminates $z_1,z_2,\ldots$ with coefficients in $\mathbb{R}$.
All identities below hold in the ring of FPS; no convergence considerations are needed.

\paragraph{Sum.}
$(f+g)(z) = \sum_n (a_n+b_n)\,z^n$.

\paragraph{Product (Cauchy convolution).}
$(f\cdot g)(z) = \sum_n \bigl(\sum_{j=0}^n a_j\,b_{n-j}\bigr) z^n$.

\paragraph{Formal derivative.}
$\PDeriv{f}{z}(z) = \sum_{n\ge 1} n\,a_n\,z^{n-1}$.

\paragraph{Composition.}
If $g(0)=0$, then $(f\circ g)(z) = f(g(z)) = \sum_n a_n\,(g(z))^n$; the chain rule holds.

\paragraph{Formal integral.}
$\int_0^z f(z')\,dz' = \sum_{n\ge 0} \frac{a_n}{n+1}\,z^{n+1}$.

\paragraph{Pointing.}
The pointing operator $z\,\PDeriv{f}{z}(z) = \sum_{n\ge 0} n\,a_n\,z^n$ marks one instance of $z$ in each monomial.

\paragraph{Standard functions.}
We use three functions throughout:
\begin{align}
\text{Exponential:}\quad E(z) &= \exp(z) = \textstyle\sum_{n\ge 0}\frac{z^n}{n!}, \label{eq:exponential_fps}\\
\text{Geometric:}\quad G(z) &= \frac{1}{1-z} = \textstyle\sum_{n\ge 0} z^n, \label{eq:geometric_fps}\\
\text{Mercator:}\quad M(z) &= -\ln(1-z) = \textstyle\sum_{n\ge 1}\frac{z^n}{n}. \label{eq:mercator_fps}
\end{align}
These satisfy $E'=E$, $G'=G^2$, $M'=G$, and $M(0)=0$.
Combinatorially, $E(z)$ generates sets (bags), $G(z)$ generates strings (sequences), and $M(z)$ generates cycles.
The identity $G(z) = E(M(z))$ encodes the cycle decomposition of permutations: every permutation decomposes uniquely into disjoint cycles.

\subsection{Canonical partition functions and factorial moments}
\label{sec:gf_canonical}

\paragraph{Canonical partition function and probabilities.}
Fix an atom composition $\mathbf{N}=(N_a)_{a\in A}$.
The canonical partition function is
\begin{equation}\label{eq:canonical_partition_general}
\GFensemble_{\mathbf{N}}
:= \Bigl(\prod_{a\in A} N_a!\Bigr)\,[\mathbf{z}^{\mathbf{N}}]\;\GFensemble,
\end{equation}
where $[\mathbf{z}^{\mathbf{N}}]$ extracts the coefficient of $\prod_a z_a^{N_a}$.
The canonical probability of a configuration $\mathbf{n}=(n_s)_{s\in\Species}$ satisfying the conservation constraints is
\begin{equation}\label{eq:canonical_probability_general}
P_{\mathbf{N}}(\mathbf{n})
:= \frac{\prod_a N_a!}{\GFensemble_{\mathbf{N}}}
\;\prod_{s\in\Species}\frac{w_s^{\,n_s}}{n_s!},
\end{equation}
where the dimensionless species weight is
\begin{equation}\label{eq:species_weight_def}
w_s = \frac{1}{\varphi(s)}\prod_{b\in B} z_b^{\beta(s,b)}.
\end{equation}
The full Poisson activity is $\lambda_s = w_s\prod_a z_a^{\alpha(s,a)}$, so $w_s$ is the $z_a$-independent part.

\paragraph{Probability generating function and shift trick.}
Introduce one marker variable $\xi_s$ per species and define
\begin{equation}\label{eq:canonical_pgf}
G_{\mathbf{N}}(\boldsymbol{\xi})
:=\sum_{\mathbf{n}}P_{\mathbf{N}}(\mathbf{n})\prod_{s\in\Species}\xi_s^{n_s}.
\end{equation}
Factorial moments are obtained by differentiation at $\boldsymbol{\xi}=\mathbf{1}$:
\[
\left.\prod_{j=1}^{m}\frac{\partial^{k_j}}{\partial \xi_{s_j}^{k_j}}
G_{\mathbf{N}}(\boldsymbol{\xi})\right|_{\boldsymbol{\xi}=\mathbf{1}}
=
\Bigl\langle \prod_{j=1}^{m}(n_{s_j})_{k_j}\Bigr\rangle_{\!\mathbf{N}}.
\]
To evaluate these derivatives, shift variables by $\boldsymbol{\xi}\mapsto\boldsymbol{\xi}+\mathbf{1}$ so derivatives are taken at $\boldsymbol{\xi}=\mathbf{0}$; this pulls down powers $w_s\,\mathbf{z}^{\boldsymbol{\alpha}(s)}$ and converts the moment computation into a depleted-composition coefficient extraction. Appendix~\ref{app:proofs} gives the detailed derivation.

\begin{theorem}[Canonical factorial moments]
\label{thm:canonical_factorial_moments}
Let $s_1,\dots,s_m\in\Species$ be species and $k_1,\dots,k_m\ge 0$.
In the canonical ensemble with atom composition~$\mathbf{N}$, the joint factorial moment of the copy numbers $n_{s_j}$ is
\begin{equation}\label{eq:canonical_factorial_moments}
\Bigl\langle \prod_{j=1}^{m}(n_{s_j})_{k_j}\Bigr\rangle_{\!\mathbf{N}}
=
\prod_{j=1}^{m}w_{s_j}^{\,k_j}
\;\cdot\;
\prod_{a\in A}(N_a)_{m_a}
\;\cdot\;
\frac{\GFensemble_{\mathbf{N}-\mathbf{m}}}{\GFensemble_{\mathbf{N}}},
\end{equation}
where $m_a = \sum_{j=1}^{m}k_j\,\alpha(s_j,a)$.
In particular, the expected copy number of species~$s$ is
\begin{equation}\label{eq:canonical_expected_value}
\langle n_s\rangle_{\mathbf{N}}
=
\frac{\prod_{b}z_b^{\beta(s,b)}}{\varphi(s)}
\;\cdot\;
\prod_{a}(N_a)_{\alpha(s,a)}
\;\cdot\;
\frac{\GFensemble_{\mathbf{N}-\boldsymbol{\alpha}(s)}}{\GFensemble_{\mathbf{N}}}.
\end{equation}
\end{theorem}

The proof, multi-index formulation, telescoping factorization, and recursive identities from partition-function ratios are collected in Appendix~\ref{app:proofs}.
Because the argument uses only the product-form structure $\GFensemble=\exp(\GFcomb)$, the theorem applies in both the acyclic and cyclic regimes.

\subsection{Saddle-point approximation and the deterministic bridge}
\label{sec:gf_saddle}

The canonical factorial moments of Theorem~\ref{thm:canonical_factorial_moments} express observables as ratios of partition functions $\GFensemble_{\mathbf{N}-\mathbf{m}}/\GFensemble_{\mathbf{N}}$.
For large atom numbers, exact coefficient extraction becomes expensive; the saddle-point method provides an asymptotic approximation whose leading order recovers the deterministic equilibrium.

\paragraph{Cauchy integral representation.}
For a single conserved atom type with total count~$N$, the canonical partition function~\eqref{eq:canonical_partition_general} is
\begin{equation}\label{eq:cauchy_integral}
\GFensemble_N
= N!\,[x^N]\,e^{\GFcomb(x)}
= \frac{N!}{2\pi i}\oint \frac{e^{\GFcomb(z)}}{z^{N+1}}\,dz,
\end{equation}
where the contour encircles the origin inside the disk of convergence of~$\GFcomb$.

\paragraph{Saddle-point equation as conservation law.}
The saddle point $x^*>0$ of the coefficient extraction $[z^N]\,e^{\GFcomb(z)}$ satisfies
\begin{equation}\label{eq:saddle_point_equation}
x^*\,\GFcomb'(x^*) = N,
\end{equation}
which is the Hayman saddle-point equation $a(x^*)=N$ for $f(z)=e^{\GFcomb(z)}$ (see Appendix~\ref{app:saddle}).
This coincides with the deterministic conservation equation~\eqref{eq:conservation_general}: the equilibrium fugacity that solves the mass-balance equation is the saddle point of the coefficient integral.

\begin{proposition}[Saddle-point formula for canonical partition functions]
\label{prop:saddle_point_general}
Let $\GFcomb(x)$ have non-negative coefficients and radius of convergence $R>0$, and let $x^*\in(0,R)$ satisfy~\eqref{eq:saddle_point_equation}.
Then the Gaussian saddle-point estimate, combined with Stirling's approximation, gives
\begin{equation}\label{eq:saddle_point_general}
\widehat{\GFensemble}_N
:=
\sqrt{\frac{N}{\lambda_N}}\;
\exp\!\bigl(\GFcomb(x^*)\bigr)\,
\left(\frac{N}{e\,x^*}\right)^{\!N},
\end{equation}
where
\begin{equation}\label{eq:lambda_def}
\lambda_N := (x^*)^2\,\GFcomb''(x^*) + N.
\end{equation}
The quantity $\lambda_N$ is the variance parameter of the Gaussian fluctuation around the saddle point; it equals the second derivative of $-h$ at~$x^*$ (in the variable $\theta$ on the circle $z=x^*e^{i\theta}$).
\end{proposition}

The derivation, including the change of variables $z=x^*e^{i\theta}$, the Gaussian evaluation, and the application of Stirling's formula, is carried out in Appendix~\ref{app:saddle}.

\paragraph{Error between exact and approximate forms.}
Expanding the integrand to fourth order around~$x^*$ and including the next Stirling correction yields the relative error
\begin{equation}\label{eq:saddle_point_error}
\frac{\GFensemble_N}{\widehat{\GFensemble}_N}
= 1
+ \frac{1}{8}\,\frac{\kappa_4}{\kappa_2^2}
- \frac{5}{24}\,\frac{\kappa_3^2}{\kappa_2^3}
+ \frac{1}{12N}
+ O\!\left(N^{-2}\right),
\end{equation}
where $\kappa_j := \frac{d^j}{du^j}h(x^*e^{u})\big|_{u=0}$ (with $u=i\theta$) are the expansion coefficients of the exponent at the saddle point, satisfying $\kappa_2 = \lambda_N$.
Since $\kappa_j = \Theta(N)$ for all $j\ge 2$, each fraction in~\eqref{eq:saddle_point_error} is $O(1/N)$, making the total relative error $O(1/N)$.
Appendix~\ref{app:saddle} derives~\eqref{eq:saddle_point_error} and gives the explicit cumulant expressions.

\paragraph{Validity conditions.}
The saddle-point approximation is rigorously justified by the \emph{H-admissibility} framework of Hayman~\cite{hayman1956generalisation}, codified by Flajolet and Sedgewick~\cite{flajolet2009analytic} (Definition~VIII.1, Theorem~VIII.4).
For assembly systems, $\GFensemble = \exp(\GFcomb)$ is H-admissible whenever $\GFcomb$ has non-negative coefficients and $b(r)\to\infty$ as $r\to R^-$, because the exponential preserves admissibility~\cite{flajolet2009analytic} (Theorem~VIII.5).
Both conditions hold for our species generating functions, whose coefficients are positive symmetry weights $1/\varphi(s)$.
The approximation requires the saddle point to lie within the radius of convergence: $x^*<R$.
In physical terms this means the system must be below the saturation or gelation threshold; Section~\ref{sec:singularity} develops tools for computing $R$ directly from the structural equations.
Appendix~\ref{app:saddle} states the H-admissibility conditions in full and verifies them for the assembly setting.

\begin{remark}[Deterministic bridge]
The saddle-point value $\GFcomb(x^*)$ is the deterministic free energy: $\log\mathcal{Z}\approx \GFcomb(x^*)$, and the Gaussian prefactor $\sqrt{N/\lambda_N}$ governs the leading finite-size correction.
The saddle-point approximation therefore bridges the stochastic (canonical) and deterministic descriptions: the deterministic equilibrium is the leading-order term, and the fluctuation parameter $\lambda_N$ quantifies departures at finite~$N$.
This bridge is illustrated in the worked examples of Sections~\ref{sec:cyclefree} and~\ref{sec:cycles}.
\end{remark}

\subsection{Singularity analysis from structural equations}
\label{sec:singularity}

The saddle-point approximation requires $x^*<R$, where $R$ is the radius of convergence of~$\GFcomb$.
At $x^*=R$ the conservation map $c(z)=\gamma\,z\,\GFcomb'(z)$ ceases to be invertible, signaling a phase boundary (gelation, saturation).
Characterizing $R$ and the nature of the singularity from the structural equations---recursions and binding equations---without solving for $\GFcomb$ in closed form is therefore the key step.

\paragraph{Smooth implicit-function schema.}
Suppose $\GFcomb$ satisfies a species recursion of the form
\begin{equation}\label{eq:implicit_recursion}
\GFcomb = F(z,\GFcomb),
\end{equation}
where $F$ is a known function of a primary variable~$z$ (typically an atom fugacity) and $\GFcomb$ itself.
Implicit differentiation gives
\begin{equation}\label{eq:implicit_derivative}
\GFcomb' = \frac{F_z}{1 - F_{\GFcomb}},
\end{equation}
which diverges when $F_{\GFcomb}=1$.
The \emph{critical system} determining the dominant singularity $z_c$ and the critical value $\tau=\GFcomb(z_c)$ is
\begin{equation}\label{eq:critical_system}
\tau = F(z_c,\tau),
\qquad
F_{\GFcomb}(z_c,\tau)=1,
\end{equation}
two equations in two unknowns that require only the recursion function~$F$, not the explicit solution for~$\GFcomb$.

\begin{proposition}[Singularity from a species recursion]
\label{prop:implicit_singularity}
Let $\GFcomb = F(z,\GFcomb)$ with $F$ analytic in both arguments, $F(0,0)=0$, $F_z(0,0)\neq 0$, and all Taylor coefficients of~$F$ non-negative.
Suppose the critical system~\eqref{eq:critical_system} has a solution $(z_c,\tau)$ with $z_c>0$.
\begin{enumerate}
\item If $F$ is nonlinear in $\GFcomb$ (i.e.\ $F_{\GFcomb\GFcomb}(z_c,\tau)\neq 0$), then $z_c$ is the radius of convergence of~$\GFcomb$ and the singularity is a square-root branch point:
\[
\GFcomb(z) \sim \tau - C\sqrt{1-z/z_c},
\qquad
[z^n]\,\GFcomb \sim C'\,z_c^{-n}\,n^{-3/2}.
\]
\item If $F$ is linear in $\GFcomb$, then $F_{\GFcomb}$ is independent of~$\GFcomb$ and the singularity is a pole (simple or higher-order) rather than a branch point.
\end{enumerate}
\end{proposition}

Case~(i) is the smooth implicit-function schema of Flajolet and Sedgewick~\cite{flajolet2009analytic} (Theorem~VII.3); case~(ii) covers linear recursions such as the linear-polymer system of Section~\ref{sec:cyclefree_linear}, where $F=a+ab\,\GFcomb$ gives $F_{\GFcomb}=ab$ and the singularity at $a=1/b$ is a simple pole.
The square-root universality of case~(i) governs all bounded-valence branching systems: for instance, the trivalent recursion $T=x(1+yT)(1+zT)$ mentioned in Section~\ref{sec:outlook} has $F_T=x(y+z+2yzT)$, and the critical system yields a finite $x_c(y,z)$ with a genuine gelation threshold.

\paragraph{Jacobian criterion from the conservation map.}
The conservation equation~\eqref{eq:conservation_general} defines a map from fugacities $(z_a)$ to total concentrations $(c_a^{\mathrm{total}})$.
Its critical points---where the map ceases to be locally invertible---are the zeros of the Jacobian:
\begin{equation}\label{eq:jacobian_criterion}
\det\!\left(\PDeriv{c_a^{\mathrm{total}}}{z_{a'}}\right) = 0.
\end{equation}
For a single atom type this reduces to
\[
\Deriv{c}{z} = \gamma\bigl(\GFcomb'(z) + z\,\GFcomb''(z)\bigr) = 0.
\]
Substituting $\GFcomb' = F_z/(1-F_{\GFcomb})$ and differentiating again expresses the Jacobian condition entirely in terms of partial derivatives of~$F$, without knowledge of~$\GFcomb$.
The critical concentration
\begin{equation}\label{eq:critical_concentration}
c_{\max} = \lim_{z\to R^-}\gamma\,z\,\GFcomb'(z)
\end{equation}
is the maximum atom count for which the saddle-point equation~\eqref{eq:saddle_point_equation} admits a solution, connecting the singularity directly to the phase boundary of Section~\ref{sec:gf_saddle}.

\paragraph{Discriminant for algebraic systems.}
When the recursion~\eqref{eq:implicit_recursion} is polynomial in $\GFcomb$---as it is for all bounded-valence assembly systems---$\GFcomb$ satisfies an algebraic equation $P(z,\GFcomb)=0$.
The singularities of $\GFcomb(z)$ are the zeros of the discriminant
\begin{equation}\label{eq:discriminant}
\Delta(z) = \operatorname{Res}_{\GFcomb}\!\left(P,\;\PDeriv{P}{\GFcomb}\right),
\end{equation}
where $\operatorname{Res}_{\GFcomb}$ denotes the resultant with respect to~$\GFcomb$.
This is a purely algebraic computation, applicable even when $\GFcomb$ cannot be expressed in closed form.
It recovers the smooth-schema conditions and extends to higher-order singularities when the recursion has degenerate structure.

\paragraph{Systems of coupled recursions.}
For multicomponent systems with coupled species generating functions $\boldsymbol{\GFcomb} = \mathbf{F}(z,\boldsymbol{\GFcomb})$, the critical condition generalizes to
\begin{equation}\label{eq:dlw_condition}
\det\!\bigl(I - \partial\mathbf{F}/\partial\boldsymbol{\GFcomb}\bigr) = 0.
\end{equation}
Under irreducibility and aperiodicity conditions on~$\mathbf{F}$ (the hypotheses of the Drmota--Lalley--Woods theorem~\cite{drmota2009random,flajolet2009analytic}), all components $\GFcomb_i$ share the same dominant singularity~$z_c$, which is again a square-root branch point with $n^{-3/2}$ coefficient asymptotics.

\begin{remark}[Cycles and singularities]
In the cyclic regime of Section~\ref{sec:cycles}, the linearized ensemble PDE~\eqref{eq:linearized_cycle_identity} propagates singularities of the initial data $\GFensemble_0$ through the bond-variable evolution.
The implicit-function schema applies to the acyclic backbone; cycle closures shift the location of the singularity (by modifying the effective recursion) but preserve the universality class.
\end{remark}

\section{Cycle-free assembly}
\label{sec:cyclefree}

This section develops the generating-function calculus for assemblies without cycles (trees, chains, branching structures).
We first present the general framework, then work out linear polymers in detail.

\subsection{General theory}
\label{sec:cyclefree_general}

We work throughout in the dimensionless variables $z_a$, $z_b$ introduced in Section~\ref{sec:deterministic}.

\paragraph{Species generating function with free-site variables.}
Fix an assembly system with atom types $A$, bond types $B$, and site types $\Sigma$.
For each species $s\in\Species$, write $\varphi(s)$ for its number of automorphisms.
The monomial associated to $s$ is
\begin{equation}\label{eq:species_monomial}
g(s) = \frac{1}{\varphi(s)}\,
\prod_{a\in A} z_a^{\alpha(s,a)}\,
\prod_{b\in B} z_b^{\beta(s,b)}\,
\prod_{\sigma\in\Sigma} z_\sigma^{\sigma_{\mathrm{free}}(s,\sigma)},
\end{equation}
and the species generating function is $\GFspecies = \sum_{s\in\Species} g(s)$, a formal power series in atom variables $z_a$, bond variables $z_b$, and auxiliary free-site variables $z_\sigma$.
Setting all site variables to~$1$ recovers the combinatorial species generating function of Section~\ref{sec:stochastic}:
\[
\GFcomb = \GFspecies\big|_{z_\sigma=1}
= \sum_{s\in\Species}\frac{1}{\varphi(s)}\,
\prod_{a\in A} z_a^{\alpha(s,a)}\,
\prod_{b\in B} z_b^{\beta(s,b)}.
\]
The weight $1/\varphi(s)$ is the generating-function expression of ``counting up to symmetry,'' grounded in the groupoid-cardinality viewpoint of Baez--Dolan \cite{baez2001finite} (see Appendix~\ref{app:foundations}).

\paragraph{General binding equation (species GF).}
Suppose bond type $b$ binds a site of type $\sigma_1$ to a site of type $\sigma_2$ (with $\sigma_1\neq\sigma_2$).
In an acyclic assembly, discounting one $b$-bond necessarily splits the assembly into two connected components.
This yields the \emph{acyclic binding equation}:
\begin{equation}\label{eq:general_binding_eq_species}
\PDeriv{\GFspecies}{z_b} = \PDeriv{\GFspecies}{z_{\sigma_1}}\,\PDeriv{\GFspecies}{z_{\sigma_2}}.
\end{equation}
The derivative on the left discounts a bond; the product on the right records the two fragments, each marked with its newly exposed free site.
When $\sigma_1=\sigma_2=\sigma$ (same-site binding), the identity acquires a symmetry factor of $1/2$.

\paragraph{System of equations for site derivatives.}
Define
\[
F_\sigma := \PDeriv{\GFspecies}{z_\sigma},\qquad \sigma\in\Sigma.
\]
In the acyclic regime, marking a free site of type $\sigma$ selects a root monomer and orients the species as a rooted tree.
For a monomer type $a\in A$, let $\Sigma(a)$ be its finite set of site positions, and let $\tau(x)\in\Sigma$ be the type of a position $x\in\Sigma(a)$.
If $x$ is the marked position ($\tau(x)=\sigma$), then each other position $x'\in\Sigma(a)\setminus\{x\}$ contributes either a free-site factor $z_{\tau(x')}$ or a bonded branch factor $z_b\,F_{\rho(b,x')}$, where $b$ is a bond type compatible with $x'$ and $\rho(b,x')$ is the partner site type reached through~$b$.
Hence
\begin{equation}\label{eq:cyclefree_Fsigma_system}
F_\sigma
=
\sum_{a\in A} z_a
\sum_{\substack{x\in\Sigma(a)\\ \tau(x)=\sigma}}
\prod_{x'\in\Sigma(a)\setminus\{x\}}
\left(
z_{\tau(x')}
+
\sum_{b\in B(x')} z_b\,F_{\rho(b,x')}
\right),
\end{equation}
with $B(x')$ the set of bond types that can bind at site position $x'$.
Thus each summand is exactly ``monomer variable times a product of factors $(\sigma'+bF_{\sigma'})$'' in the sense that every factor is ``free-site term + bonded-subtree term.''

When each site type has a unique partner map $\pi:\Sigma\to\Sigma$ and one associated bond variable $z_{b_{\sigma'}}$, \eqref{eq:cyclefree_Fsigma_system} simplifies to
\begin{equation}\label{eq:cyclefree_Fsigma_system_unique_partner}
F_\sigma
=
\sum_{a\in A} z_a
\sum_{\substack{x\in\Sigma(a)\\ \tau(x)=\sigma}}
\prod_{x'\in\Sigma(a)\setminus\{x\}}
\left(z_{\tau(x')} + z_{b_{\tau(x')}}\,F_{\pi(\tau(x'))}\right).
\end{equation}
In same-type binding models ($\pi(\sigma')=\sigma'$), each factor has the literal form $z_{\sigma'} + z_{b_{\sigma'}}F_{\sigma'}$.

Setting all free-site variables to $1$ gives a closed polynomial system for reduced derivatives
\[
\Psi_\sigma := F_\sigma\big|_{z_\eta=1\;\forall\eta\in\Sigma},
\]
namely
\begin{equation}\label{eq:cyclefree_Psi_system}
\Psi_\sigma
=
\sum_{a\in A} z_a
\sum_{\substack{x\in\Sigma(a)\\ \tau(x)=\sigma}}
\prod_{x'\in\Sigma(a)\setminus\{x\}}
\left(
1
+
\sum_{b\in B(x')} z_b\,\Psi_{\rho(b,x')}
\right).
\end{equation}
This is the general mechanism behind rooted-species recursions in cycle-free systems and the coupled-recursion singularity criteria of Section~\ref{sec:singularity}.

\paragraph{Ensemble generating function.}
The ensemble generating function is
\begin{equation}\label{eq:ensemble_gf_def}
\GFensemble = \exp(\GFcomb).
\end{equation}
Exponentiation forms multisets (``ensembles'') of species with factorial factors for indistinguishability.
No system-size parameter $\gamma$ appears: all variables are dimensionless.

\paragraph{Binding equation (ensemble GF).}
By the chain rule applied to \eqref{eq:ensemble_gf_def} and \eqref{eq:general_binding_eq_species}:
\begin{equation}\label{eq:general_binding_eq_ensemble}
\PDeriv{\GFensemble}{z_b}
= \GFensemble \cdot \left.\PDeriv{\GFspecies}{z_{\sigma_1}}\,\PDeriv{\GFspecies}{z_{\sigma_2}}\right|_{z_\sigma=1}.
\end{equation}
This equation \emph{couples} the ensemble GF to the species GF.
It is linear in $\GFensemble$ but does \emph{not} close as a PDE in $\GFensemble$ alone.
Writing $\widetilde{\GFensemble} = \exp(\GFspecies)$ for the ensemble GF retaining site variables, the mixed partial satisfies
\begin{equation}\label{eq:general_binding_eq_ensemble_corrected}
\PDeriv{\GFensemble}{z_b}
= \left.\frac{\partial^2\widetilde{\GFensemble}}{\partial z_{\sigma_1}\,\partial z_{\sigma_2}}\right|_{z_\sigma=1}
- \GFensemble\,\left.\frac{\partial^2\GFspecies}{\partial z_{\sigma_1}\,\partial z_{\sigma_2}}\right|_{z_\sigma=1},
\end{equation}
and the residual term prevents closure.
Linearization requires the cycle-opening term of Section~\ref{sec:cycles}.

\paragraph{Species recursion for rooted structures.}
When assemblies have a distinguished root, the species GF may admit a recursion.
In a system with atom $z_x$ and bond $z_y$, a rooted linear polymer is either a single atom or an atom bonded to another rooted polymer: $T = z_x + z_x z_y T$.
In branching systems with multiple bond types, the recursion involves products of factors $(1+\text{bond}\cdot T)$.

\paragraph{Ensemble recursion by integration.}
Integrating the binding equation~\eqref{eq:general_binding_eq_ensemble} with respect to $z_b$:
\begin{equation}\label{eq:ensemble_recursion_general}
\GFensemble(z_b) = \GFensemble(0) + \int_0^{z_b} \GFensemble\,\left.\PDeriv{\GFspecies}{z_{\sigma_1}}\,\PDeriv{\GFspecies}{z_{\sigma_2}}\right|_{z_\sigma=1}\,dz_b'.
\end{equation}
Because $\GFcomb$ is independently determined, this is a linear integral equation in $\GFensemble$.

\paragraph{Equilibrium concentrations.}
The generating-function framework provides a systematic route from binding parameters to molar concentrations.

\emph{Pointing and conservation.}
The \emph{pointing operator} $z_a\,\partial/\partial z_a$ marks one atom of type~$a$ in each species.
The total molar concentration of atoms of type~$a$ is
\begin{equation}\label{eq:conservation_general}
c_a^{\mathrm{total}} = \gamma\,z_a\,\PDeriv{\GFcomb}{z_a},
\end{equation}
where $\gamma = 1/(N_A V)$.
This is the conservation equation: it relates the known total concentration to the unknown dimensionless free-monomer variable~$z_a$.
Given the bond variables $z_b$, one solves~\eqref{eq:conservation_general} for each~$z_a$.

\emph{Species concentrations.}
Once the free-monomer variables $z_a$ are determined, the molar concentration of species~$s$ is
\begin{equation}\label{eq:concentration_general}
c_s = \gamma\,\frac{1}{\varphi(s)}\,
\prod_{a\in A} z_a^{\alpha(s,a)}\,
\prod_{b\in B} z_b^{\beta(s,b)}.
\end{equation}
The universal prefactor $\gamma$ follows from the Euler characteristic identity (Section~\ref{sec:euler-characteristic}), which ensures that every connected species---acyclic or cyclic---carries exactly one net power of $\gamma$ once atom and bond dimensions are accounted for.
In the molar-variable representation the concentration formula reads $c_s = \gamma^{\vartheta(s)}\,\varphi(s)^{-1}\prod_a c_a^{\alpha(s,a)}\prod_b y_b^{\beta(s,b)}$, where the topology-dependent factor $\gamma^{\vartheta(s)}$ is explicit; for acyclic species $\vartheta(s)=0$ and this factor is~$1$.

\emph{Change of variables from molar inputs.}
In practice, binding constants are reported as molar dissociation constants $K_{d,b}$ (units of concentration) and total amounts as molar concentrations $c_a^{\mathrm{total}}$.
The conversion to dimensionless variables is
\begin{equation}\label{eq:molar_to_dimensionless}
z_b = \frac{\gamma}{K_{d,b}} = \gamma\,y_b,
\qquad
z_a = \frac{c_a}{\gamma},
\end{equation}
where $y_b = 1/K_{d,b}$ is the molar binding affinity and $c_a$ is the free-monomer molar concentration.
Substituting into~\eqref{eq:conservation_general}, one solves for $z_a$ (equivalently $c_a = \gamma\,z_a$), then evaluates~\eqref{eq:concentration_general} for each species.

\emph{Converting molar affinities to dimensionless.}
Computing $\GFensemble_{\mathbf{N}}$ from experimental binding constants requires converting to dimensionless variables via~\eqref{eq:molar_to_dimensionless}.
Given molar binding affinities $y_b = \kappa_b^+/\kappa_b^-$, the dimensionless bond variable is $z_b = \gamma\,y_b$.
In the stochastic description, the bimolecular rate constant becomes $\kappa_b^+/(N_A V) = \gamma\,\kappa_b^+$, so the dimensionless affinity $z_b = \gamma\,\kappa_b^+/\kappa_b^-$ is exactly the ratio of stochastic rate constants.
Canonical partition functions and factorial moments are developed in the shared-theory discussion of Section~\ref{sec:gf_canonical}; here we now specialize to acyclic binding identities and examples.

\subsection{Linear polymers: cycle-free}
\label{sec:cyclefree_linear}

\paragraph{Bond system.}
The linear-polymer system has one atom variable $a$, one bond variable $b$, and two site variables $p,q$.
Each monomer contributes one $p$-site and one $q$-site; each bond consumes one $p$-site and one $q$-site.
Valid assemblies are finite paths (no cycles, no branching).

\paragraph{Species generating function.}
A polymer of length $n$ has $n$ atoms, $n-1$ bonds, one free $p$-site and one free $q$-site.
All polymers have trivial automorphisms ($\varphi=1$).
The species generating function is
\[
\GFspecies = pq\sum_{n=1}^\infty a^n b^{n-1} = \frac{pqa}{1-ab}.
\]

\paragraph{Recursion.}
The species GF satisfies:
\begin{equation}\label{eq:linear_recursion}
\GFspecies = pqa + ab\GFspecies.
\end{equation}
A polymer is either a single monomer (with both sites free) or a monomer bonded to a polymer.
Setting $p=q=1$ gives
\[
\GFcomb = \frac{a}{1-ab}.
\]

\paragraph{Binding equation.}
The acyclic binding equation \eqref{eq:general_binding_eq_species} reads
\begin{equation}\label{eq:linear_polymer_discounting_identity}
\PDeriv{\GFspecies}{b} = \PDeriv{\GFspecies}{p}\,\PDeriv{\GFspecies}{q},
\end{equation}
since discounting a bond splits the polymer into two fragments.
Direct computation confirms:
\[
\PDeriv{\GFspecies}{b}
= \frac{pqa^2}{(1-ab)^2}
= \PDeriv{\GFspecies}{p}\,\PDeriv{\GFspecies}{q}.
\]

The formal derivative $\partial/\partial z$ acts by \emph{discounting} one instance of~$z$: for bond variables, this is fission (removing a bond splits a connected assembly); for atom variables, it marks (``points to'') one atom.

\paragraph{Equilibrium concentrations.}
The dimensionless monomial for a polymer of length $n$ is $a^n b^{n-1}$.
By \eqref{eq:concentration_general}, its molar concentration is
\[
c_n = \gamma\,a^n b^{n-1}.
\]
With one atom type, the conservation equation \eqref{eq:conservation_general} becomes
\begin{equation}\label{eq:linear_conservation}
c^{\mathrm{total}} = \gamma\,a\,\PDeriv{\GFcomb}{a} = \gamma\,\frac{a}{(1-ab)^2}.
\end{equation}
This is the \emph{pointing} operator $a\,\partial_a$ applied to $\GFcomb$: it marks one atom in each species, equivalent to rooting a polymer.
In molar variables (\eqref{eq:molar_to_dimensionless}), with $c_{\mathrm{free}}=\gamma a$ and $y=1/K_d$ so that $b=\gamma y$, this is
\[
c^{\mathrm{total}} = \frac{c_{\mathrm{free}}}{(1-c_{\mathrm{free}}y)^2},
\]
and therefore
\[
c_{\mathrm{free}}
= \frac{1}{c^{\mathrm{total}}}
\left(\frac{1-\sqrt{1+4c^{\mathrm{total}}y}}{2y}\right)^2.
\]

\paragraph{Alternative recursion via integration.}
Integrating the binding equation:
\begin{equation}\label{eq:linear_polymer_integral_recursion}
\GFspecies = pqa + \int_0^b \PDeriv{\GFspecies}{p}\,\PDeriv{\GFspecies}{q}\,db'.
\end{equation}
The formal integral adds one instance of $b$ and divides by the new exponent, ensuring the added bond is undistinguished.
This reads: ``a polymer is either a single monomer, or two fragments---one with a free $p$-site and one with a free $q$-site---joined by a bond.''

\paragraph{Ensemble generating function.}
Setting $p=q=1$ gives $\GFcomb = a/(1-ab)$ and therefore
\[
\GFensemble = \exp\!\left(\frac{a}{1-ab}\right).
\]
By the chain rule:
\begin{equation}\label{eq:linear_ensemble_binding}
\PDeriv{\GFensemble}{b}
= \GFensemble\,\left.\PDeriv{\GFspecies}{p}\,\PDeriv{\GFspecies}{q}\right|_{p=q=1}
= \GFensemble\cdot\GFcomb^2.
\end{equation}
Integrating:
\begin{equation}\label{eq:linear_ensemble_recursion}
\GFensemble = e^{a} + \int_0^b \GFensemble\,\GFcomb^2\,db'.
\end{equation}
An ensemble of polymers is either an ensemble of free atoms, or an ensemble with two polymers joined.

\paragraph{Canonical partition function.}
The canonical partition function for $N$ total atoms is
\begin{equation}\label{eq:linear_canonical_partition}
\GFensemble_N = N!\,[a^N]\;\exp\!\left(\frac{a}{1-ab}\right).
\end{equation}
By the Fa\`a di Bruno formula for the exponential of an EGF, with $\hat t_n = n!\,b^{n-1}$:
\begin{equation}\label{eq:linear_faa_di_bruno}
\GFensemble_N
= \sum_{\pi\in\Pi([N])}\;\prod_{B\in\pi} |B|!\,b^{|B|-1},
\end{equation}
where $\Pi([N])$ is the lattice of set partitions of $[N]$ and $|\pi|$ is the number of blocks.
A closed form is
\begin{equation}\label{eq:linear_explicit}
\GFensemble_N
= \sum_{i=0}^{N-1} \frac{N!}{(N-i)!}\,\binom{N-1}{i}\,b^{\,i}.
\end{equation}

\paragraph{Canonical probabilities.}
In the canonical ensemble with $N$ atoms, the probability of an ensemble $n=(n_i)_{i\ge 1}$ with $\sum_i i\,n_i = N$ is
\[
P(n) = \frac{N!}{\GFensemble_N}\prod_{i\ge 1}\frac{(b^{i-1})^{n_i}}{n_i!}.
\]

\section{Assembly with cycles}
\label{sec:cycles}

This section extends the calculus to assemblies containing cycles (rings).
The species generating function retains the same $\gamma$-free form as in the acyclic case: because the Euler characteristic identity guarantees that every connected species carries exactly one net power of~$\gamma$ in molar variables, the dimensionless monomials $g(s)$ of \eqref{eq:species_monomial} absorb the cycle rank automatically.
What changes is the \emph{binding equation}: discounting a bond in a cyclic species may open a cycle rather than split the assembly, introducing an additional cycle-opening term.

\subsection{General theory}
\label{sec:cycles_general}

\paragraph{Species generating function.}
The species GF has the same form as in the acyclic case:
\begin{equation}\label{eq:cycle_weighted_species_gf}
\GFspecies = \sum_{s\in\Species} g(s),
\end{equation}
where $g(s)$ is the symmetry-weighted monomial \eqref{eq:species_monomial}.
No cycle-rank weighting appears: the dimensionless variables $z_a$, $z_b$ absorb all factors of~$\gamma$.

\paragraph{Cyclic discounting identity.}
Fix a bond type $z_b$ binding site $z_{\sigma_1}$ to $z_{\sigma_2}$.
Discounting a $z_b$-bond in an acyclic species produces fission.
In a cyclic species there is a second possibility: discounting may \emph{open a cycle} while keeping the structure connected.
This yields:
\begin{equation}\label{eq:cyclic_discounting_identity}
\PDeriv{\GFspecies}{z_b}
=
\PDeriv{\GFspecies}{z_{\sigma_1}}\,\PDeriv{\GFspecies}{z_{\sigma_2}}
+
\frac{\partial^2\GFspecies}{\partial z_{\sigma_1}\,\partial z_{\sigma_2}}.
\end{equation}
The first term is fission (acyclic contribution); the second is cycle opening.
For same-site binding ($z_{\sigma_1}=z_{\sigma_2}=z_\sigma$):
\begin{equation}\label{eq:cyclic_discounting_samesite}
\PDeriv{\GFspecies}{z_b}
=
\frac{1}{2}\left(\PDeriv{\GFspecies}{z_\sigma}\right)^{\!2}
+
\frac{1}{2}\,\frac{\partial^2\GFspecies}{\partial z_\sigma^2}.
\end{equation}

\paragraph{Linearization via the ensemble generating function.}
Equation~\eqref{eq:cyclic_discounting_identity} is nonlinear in $\GFspecies$.
Define $\widetilde{\GFensemble} = \exp(\GFspecies)$, the site-augmented ensemble GF introduced in~\eqref{eq:general_binding_eq_ensemble_corrected}.
Differentiating and using \eqref{eq:cyclic_discounting_identity} yields a \emph{linear} identity:
\begin{equation}\label{eq:linearized_cycle_identity}
\PDeriv{\widetilde{\GFensemble}}{z_b}
=
\frac{\partial^2\widetilde{\GFensemble}}{\partial z_{\sigma_1}\,\partial z_{\sigma_2}}.
\end{equation}
This is a closed linear PDE in $\widetilde{\GFensemble}$---the hallmark of the cyclic regime.
The cycle-opening term in the species equation exactly cancels the residual in \eqref{eq:general_binding_eq_ensemble_corrected} that prevented closure in the acyclic case.

\begin{remark}
This linearization is the paper's central structural result: the nonlinear species-level binding equation (a Riccati-type PDE) is linearized by the exponential transform $\GFspecies\mapsto\widetilde{\GFensemble} = \exp(\GFspecies)$, but \emph{only} when the cycle-opening term is present.
In the acyclic case, the same transform produces a residual \eqref{eq:general_binding_eq_ensemble_corrected} that prevents closure.
\end{remark}

\paragraph{Operator-exponential solution.}
Writing $\widetilde{\GFensemble}_0$ for the site-augmented ensemble GF with no $z_b$-bonds:
\begin{equation}\label{eq:operator_exponential_solution}
\widetilde{\GFensemble} = \exp\!\left(z_b\,\frac{\partial^2}{\partial z_{\sigma_1}\,\partial z_{\sigma_2}}\right)\widetilde{\GFensemble}_0.
\end{equation}
Expanding:
\begin{equation}\label{eq:operator_exponential_series}
\widetilde{\GFensemble}
=
\sum_{n=0}^\infty \frac{z_b^{\,n}}{n!}\,D^n\,\widetilde{\GFensemble}_0,
\qquad D = \frac{\partial^2}{\partial z_{\sigma_1}\,\partial z_{\sigma_2}}.
\end{equation}
The $n$-th term discounts $n$ free sites of each type and adds $n$ indistinguishable $z_b$-bonds, providing a controlled truncation scheme.

\paragraph{Ensemble recursion.}
Integrating \eqref{eq:linearized_cycle_identity}:
\begin{equation}\label{eq:ensemble_recursion_cyclic}
\widetilde{\GFensemble} = \widetilde{\GFensemble}_0 + \int_0^{z_b} \frac{\partial^2\widetilde{\GFensemble}}{\partial z_{\sigma_1}\,\partial z_{\sigma_2}}\,dz_b'.
\end{equation}

\paragraph{Multiple bond types and the aggregate bond operator.}
In systems with many bond types, each has the same fission-versus-cycle-opening structure,
and the exponential transform linearizes each one.
For heterotypic bonds ($\sigma_1\neq\sigma_2$), the linearized identity~\eqref{eq:linearized_cycle_identity} reads $\partial\widetilde{\GFensemble}/\partial z_b = D_b\,\widetilde{\GFensemble}$ with $D_b = \partial^2/(\partial z_{\sigma_1}\,\partial z_{\sigma_2})$.
For homotypic bonds ($\sigma_1=\sigma_2=\sigma$), the same-site identity~\eqref{eq:cyclic_discounting_samesite} linearizes to $\partial\widetilde{\GFensemble}/\partial z_b = \tfrac{1}{2}\,D_b\,\widetilde{\GFensemble}$ with $D_b = \partial^2/\partial z_\sigma^2$.
Because each $D_b$ is a constant-coefficient differential operator in the site variables, the operators commute:
$[D_{b_1},D_{b_2}]=0$ for all $b_1,b_2\in B$, regardless of whether bond types share sites.
Writing $c_b = 1$ for heterotypic bonds and $c_b = 1/2$ for homotypic bonds, the joint solution is
\begin{equation}\label{eq:aggregate_operator_exponential}
\widetilde{\GFensemble}
= \exp\!\left(\sum_{b\in B} c_b\,z_b\,D_b\right)\widetilde{\GFensemble}_0,
\end{equation}
and the product of operator exponentials $\prod_b\exp(c_b\,z_b D_b)$ equals the exponential of the sum.

\paragraph{Simple form of the initial data.}
With all bond variables set to zero, every species is a single monomer.
For atom type $a\in A$, let $m_\sigma(a) = |\{x\in\Sigma(a):\tau(x)=\sigma\}|$ be the number of site positions of type~$\sigma$, and define the \emph{monomer monomial}
\begin{equation}\label{eq:monomer_monomial}
\lambda_a := z_a\prod_{\sigma\in\Sigma} z_\sigma^{m_\sigma(a)}.
\end{equation}
The bond-free species generating function is $\GFspecies_0 = \sum_{a\in A}\lambda_a$, and the initial data is
\begin{equation}\label{eq:Z0_explicit}
\widetilde{\GFensemble}_0 = \exp\!\left(\sum_{a\in A}\lambda_a\right).
\end{equation}
The exponent is a finite sum of monomials.

\paragraph{Multinomial expansion.}
By commutativity, \eqref{eq:aggregate_operator_exponential} expands as
\begin{equation}\label{eq:multinomial_operator_expansion}
\widetilde{\GFensemble}
= \sum_{\mathbf{n}\ge 0}
\frac{\mathbf{c}^{\mathbf{n}}\,\mathbf{z}_b^{\,\mathbf{n}}}{\mathbf{n}!}\,\mathbf{D}^{\mathbf{n}}\,\widetilde{\GFensemble}_0,
\end{equation}
where $\mathbf{c}^{\mathbf{n}} = \prod_b c_b^{n_b}$, $\mathbf{z}_b^{\mathbf{n}} = \prod_b z_b^{n_b}$, $\mathbf{n}! = \prod_b n_b!$,
and $\mathbf{D}^{\mathbf{n}} = \prod_b D_b^{n_b}$ (in any order).
Setting $z_\sigma=1$ recovers the combinatorial ensemble generating function~$\GFensemble$.

\begin{theorem}[Canonical partition functions in the cyclic regime]
\label{thm:canonical_cyclic}
Let $\mathbf{N}=(N_a)_{a\in A}$ be an atom composition,
$M_\sigma = \sum_{a\in A} N_a\,m_\sigma(a)$ the total number of free sites of type~$\sigma$ in a pool of $\mathbf{N}$ monomers,
$D_b$ the operator from the linearized identity, and $c_b$ the coefficient from~\eqref{eq:aggregate_operator_exponential} ($c_b=1$ heterotypic, $c_b=1/2$ homotypic).
Then the canonical partition function~\eqref{eq:canonical_partition_general} is
\begin{equation}\label{eq:canonical_cyclic_operator}
\GFensemble_{\mathbf{N}}
= \left.\exp\!\left(\sum_{b\in B} c_b\,z_b\,D_b\right)
\prod_{\sigma\in\Sigma} z_\sigma^{M_\sigma}\right|_{z_\sigma=1}.
\end{equation}
For heterotypic bonds ($\sigma_1(b)\neq\sigma_2(b)$ for all~$b$, so $c_b=1$), this expands to
\begin{equation}\label{eq:canonical_cyclic_explicit}
\GFensemble_{\mathbf{N}}
= \sum_{\mathbf{n}\ge 0}
\frac{\mathbf{z}_b^{\,\mathbf{n}}}{\mathbf{n}!}
\prod_{\sigma\in\Sigma}(M_\sigma)_{\Delta_\sigma(\mathbf{n})}\,,
\end{equation}
where $(M)_k = M!/(M-k)!$ is the falling factorial,
\[
\Delta_\sigma(\mathbf{n}) = \sum_{\substack{b\in B\\ \sigma\in\{\sigma_1(b),\,\sigma_2(b)\}}} n_b
\]
counts the bond-ends consuming a site of type~$\sigma$,
and the sum truncates at $\Delta_\sigma(\mathbf{n})\le M_\sigma$ for every~$\sigma$.
\end{theorem}

\begin{proof}
The operators $D_b$ act on site variables and do not change atom-variable content,
so the atom-coefficient extraction $[\mathbf{z}_a^{\mathbf{N}}]$ commutes with $\exp(\sum_b c_b\,z_b D_b)$.
From \eqref{eq:Z0_explicit},
\[
[\mathbf{z}_a^{\mathbf{N}}]\;\widetilde{\GFensemble}_0
= \prod_{a\in A}\frac{(\prod_\sigma z_\sigma^{m_\sigma(a)})^{N_a}}{N_a!}
= \frac{1}{\prod_a N_a!}\prod_{\sigma\in\Sigma} z_\sigma^{M_\sigma}.
\]
Therefore
\[
\GFensemble_{\mathbf{N}}
= \Bigl(\prod_a N_a!\Bigr)\,[\mathbf{z}_a^{\mathbf{N}}]\;\GFensemble
= \left.\exp\!\left(\sum_b c_b\,z_b D_b\right)\prod_\sigma z_\sigma^{M_\sigma}\right|_{z_\sigma=1}.
\]
For the expansion, each $D_b^{n_b}$ differentiates $z_{\sigma_1}^{M_{\sigma_1}}$ and $z_{\sigma_2}^{M_{\sigma_2}}$ independently,
producing the falling factorial $(M_{\sigma_1} - \text{prior uses})_{n_b}$ for each site type.
Since the operators commute, the total depletion of site type~$\sigma$ is $\Delta_\sigma(\mathbf{n})$,
and evaluating the remaining powers at $z_\sigma=1$ gives~\eqref{eq:canonical_cyclic_explicit}.
\end{proof}

\begin{remark}
The falling-factorial structure has a transparent combinatorial reading.
For each bond type~$b$, choose $n_b$ sites of type $\sigma_1(b)$ and $n_b$ of type $\sigma_2(b)$ from the monomer pool,
then form a perfect matching between the two chosen sets ($n_b!$ matchings, cancelled by the $1/n_b!$ for indistinguishable bonds).
The product of falling factorials enforces that sites are consumed without replacement across all bond types.
Expression~\eqref{eq:canonical_cyclic_explicit} is a finite polynomial in the bond variables---it terminates when any site type is exhausted---and provides a direct route to canonical partition functions without first solving for the species generating function or performing coefficient extraction from the grand-canonical series.
\end{remark}

\subsection{Linear polymers with rings}
\label{sec:cycles_linear}

We revisit linear polymers, now allowing ring closure.
The same bond system $\mathsf{LinearPolymers}$ is used (atom~$a$, bond~$b$, sites~$p,q$), but validity is relaxed to allow both paths and simple cycles.

\paragraph{Species generating function.}
The species GF has two parts: $\GFspecies = \GFspecies_l + \GFspecies_r$, where $\GFspecies_l$ is the linear part and $\GFspecies_r$ the ring part.
A ring of length $n$ has $n$ atoms, $n$ bonds, $n$ rotational symmetries ($\varphi=n$), cycle rank $\vartheta=1$, and no free sites.
Thus:
\[
\GFspecies_r = \sum_{n\ge 1}\frac{(ab)^n}{n} = -\ln(1-ab) = M(ab).
\]
The appearance of the Mercator series $M(ab)$ reflects the rotational symmetry of cycles.
Collecting:
\begin{equation}\label{eq:polymer_with_rings}
\GFspecies = \frac{pqa}{1-ab} - \ln(1-ab).
\end{equation}
Setting $p=q=1$ gives the combinatorial species GF:
\[
\GFcomb = \frac{a}{1-ab} - \ln(1-ab).
\]

\paragraph{Derivatives and the species recursion.}
The ring part satisfies $\PDeriv{\GFspecies_r}{b} = \GFcomb_l$, where $\GFcomb_l = a/(1-ab)$ is the linear combinatorial GF: removing a bond from a ring yields a linear polymer.
Conversely, integrating $\GFcomb_l$ with respect to $b$ produces $\GFspecies_r$: adding a bond to a linear polymer and dividing by the new symmetry gives a ring.
The key relation is
\[
\PDeriv{\GFcomb}{b} = a + a^2\PDeriv{\GFcomb}{a},
\]
which yields the recursion:
\begin{equation}\label{eq:ring_species_recursion}
\GFcomb = a + ab + \int_0^b a^2\PDeriv{\GFcomb}{a}\,db'.
\end{equation}

\paragraph{Cyclic binding equation.}
The cyclic discounting identity \eqref{eq:cyclic_discounting_identity} for this system reads:
\[
\PDeriv{\GFspecies}{b} = \PDeriv{\GFspecies}{p}\,\PDeriv{\GFspecies}{q} + \frac{\partial^2 \GFspecies}{\partial p\,\partial q}.
\]
Setting site variables to~1: the fission term gives $\GFcomb_l^2$ and the cycle-opening term gives $\GFcomb_l$.

\paragraph{Equilibrium concentrations.}
The conservation equation \eqref{eq:conservation_general} becomes
\[
c^{\mathrm{total}} = \gamma\,a\,\PDeriv{\GFcomb}{a} = \gamma\left[\frac{a}{(1-ab)^2} + \frac{ab}{1-ab}\right].
\]
Solving for $a$ with $N = c^{\mathrm{total}}/\gamma$:
\[
a = \frac{(2N+1)b + 1 - \sqrt{(1+b)^2 + 4Nb}}{2(N+1)b^2}.
\]
The equilibrium concentrations of linear polymers and rings of length $n$ are:
\[
c_n = \gamma\,a^n b^{n-1}, \quad r_n = \gamma\,\frac{a^n b^n}{n}.
\]

\paragraph{Ensemble generating function and linearized PDE.}
The ensemble GF is $\GFensemble=\exp(\GFcomb)$:
\[
\GFensemble = \frac{1}{1-ab}\exp\!\left(\frac{a}{1-ab}\right).
\]
The factor $1/(1-ab)$ arises because $\exp(-\ln(1-ab)) = 1/(1-ab)$.
The linearized identity \eqref{eq:linearized_cycle_identity} becomes
\begin{equation}\label{eq:rings_linearized}
\PDeriv{\widetilde{\GFensemble}}{b} = \frac{\partial^2 \widetilde{\GFensemble}}{\partial p\,\partial q},
\end{equation}
a closed linear PDE in $\widetilde{\GFensemble} = \exp(\GFspecies)$.

\paragraph{Ensemble recursion.}
From the species-level relation and the chain rule:
\[
\PDeriv{\GFensemble}{b} = a\GFensemble + a^2\PDeriv{\GFensemble}{a}.
\]
Integrating:
\begin{equation}\label{eq:rings_ensemble_recursion}
\GFensemble = e^{a} + \int_0^b \left(a\GFensemble + a^2\PDeriv{\GFensemble}{a}\right)db'.
\end{equation}

\paragraph{Canonical partition function.}
\begin{equation}\label{eq:rings_partition}
\GFensemble_N = N!\,[a^N]\;\frac{1}{1-ab}\exp\!\left(\frac{a}{1-ab}\right).
\end{equation}
Expanding and extracting $[a^N]$:
\begin{equation}\label{eq:rings_explicit}
\GFensemble_N
=
\sum_{i=0}^{N}\frac{N!}{(N-i)!}\,\binom{N}{i}\,b^i.
\end{equation}
Comparing with the acyclic case \eqref{eq:linear_explicit}, the only changes are that the sum extends to $i=N$ (rather than $N-1$) and $\binom{N}{i}$ replaces $\binom{N-1}{i}$.
This is a special case of Theorem~\ref{thm:canonical_cyclic}: one atom type with $m_p(a)=m_q(a)=1$ gives $M_p=M_q=N$, and \eqref{eq:canonical_cyclic_explicit} becomes $\GFensemble_N = \sum_{n=0}^{N}\binom{N}{n}^2\,n!\;b^n$, which equals~\eqref{eq:rings_explicit} since $\binom{N}{n}^2 n! = \frac{N!}{(N-n)!}\binom{N}{n}$.

\paragraph{Canonical probabilities.}
A state is described by sequences $\mathbf{n}^{(l)}$ (linear polymers) and $\mathbf{n}^{(r)}$ (rings).
The probability is:
\[
P(\mathbf{n}^{(l)},\mathbf{n}^{(r)})
= \frac{N!}{\GFensemble_N}\,
\prod_{i\ge 1}\frac{(b^{i-1})^{\,n_i^{(l)}}}{n_i^{(l)}!}\,
\prod_{j\ge 1}\frac{(b^j/j)^{\,n_j^{(r)}}}{n_j^{(r)}!}\,.
\]

\subsection{Ring closure variants}
\label{sec:ring_closure_variants}

\paragraph{Modified ring-closure parameter.}
In some models, ring closure uses an effective dimensionless parameter $\zeta$ different from the standard value $\zeta=1$:
\begin{equation}\label{eq:cycles_with_c_solution}
\GFcomb = \frac{a}{1-ab} - \zeta\ln(1-ab).
\end{equation}

\paragraph{Separating acyclic backbone from cyclic closures.}
A practical device is to compute first an acyclic generating function $L$ and then apply a linear closure operator.
Integrating with respect to $\zeta$:
\[
\GFspecies = L + \int_0^{\zeta} \dd \zeta'\,\int_0^b \dd b'\,\frac{\partial^2\GFspecies}{\partial p\,\partial q},
\]
which unfolds as a formal operator series:
\[
\GFspecies = \sum_{n=0}^\infty \left(\int_0^{\zeta} \dd \zeta'\,\int_0^b \dd b'\,\frac{\partial^2}{\partial p\,\partial q}\right)^n L.
\]
The first few coefficients are recorded in Appendix~\ref{app:closure-expansions}.

\section{Cross-linked polymers: worked example}
\label{sec:crosslinked}

This section works through both the cycle-free and cyclic regimes for a bond system that exhibits genuine gelation: linear chains joined by cross-links.
The system is the simplest model that combines heterotypic chain-forming bonds with homotypic cross-links, and it illustrates the full generating-function toolkit developed in Sections~\ref{sec:cyclefree} and~\ref{sec:cycles}.

\subsection{Bond system}
\label{sec:crosslinked_bond_system}

The cross-linked polymer system has one atom type~$a$, three site types $p,q,r$, and two bond types:
\begin{itemize}
\item $b$ (heterotypic): binds a $p$-site to a $q$-site (chain bonds),
\item $d$ (homotypic): binds an $r$-site to an $r$-site (cross-links).
\end{itemize}
Each monomer carries one site of each type, so its maximum degree is three: one $b$-bond at~$p$, one at~$q$, and one $d$-bond at~$r$.
The $b$-bonds alone produce the linear polymers of Sections~\ref{sec:cyclefree_linear} and~\ref{sec:cycles_linear}; the $d$-bonds cross-link these chains into branching, tree-like (or cyclic) structures.

\subsection{Cycle-free case}
\label{sec:crosslinked_acyclic}

\paragraph{Site derivative system.}
Applying the general system~\eqref{eq:cyclefree_Psi_system} to this bond system, each site derivative marks a free site on the root monomer.
The remaining sites contribute factors ``free $+$ bonded subtree'':
\begin{align}
\Psi_p &= a(1+b\,\Psi_p)(1+d\,\Psi_r), \label{eq:xl_Psi_p}\\
\Psi_q &= a(1+b\,\Psi_q)(1+d\,\Psi_r), \label{eq:xl_Psi_q}\\
\Psi_r &= a(1+b\,\Psi_q)(1+b\,\Psi_p). \label{eq:xl_Psi_r}
\end{align}
In~\eqref{eq:xl_Psi_p}, the factor $(1+b\,\Psi_p)$ arises because the $q$-site is either free or bonds via~$b$ to a subtree rooted at a free $p$-site; the factor $(1+d\,\Psi_r)$ arises because the $r$-site is either free or bonds via~$d$ to a subtree rooted at a free $r$-site.
The equations for $\Psi_p$ and $\Psi_q$ are identical, so $\Psi_p=\Psi_q$.
Write $\Psi$ for this common value.
The reduced system is
\begin{equation}\label{eq:xl_reduced_system}
\Psi = a(1+b\,\Psi)(1+d\,\Psi_r),
\qquad
\Psi_r = a(1+b\,\Psi)^2.
\end{equation}

\paragraph{Implicit cubic.}
Substituting $u = 1+b\,\Psi$ (so $\Psi=(u-1)/b$ and $\Psi_r = a\,u^2$) and eliminating $\Psi_r$:
\begin{equation}\label{eq:xl_cubic}
u(1-ab) - a^2 b d\,u^3 = 1.
\end{equation}
At $d=0$ this gives $u=1/(1-ab)$, recovering the linear-polymer site derivative $\Psi=a/(1-ab)$ of Section~\ref{sec:cyclefree_linear}.

\paragraph{Species generating function by integration.}
The acyclic binding equations~\eqref{eq:general_binding_eq_species} specialize to
\begin{equation}\label{eq:xl_binding_eqs}
\PDeriv{\GFcomb}{b} = \Psi^2,
\qquad
\PDeriv{\GFcomb}{d} = \frac{1}{2}\,\Psi_r^2,
\end{equation}
where the factor $1/2$ in the second equation reflects the symmetry of same-site binding.
Integration proceeds along a path in $(b,d)$-space.
\emph{Step~1}: at $d=0$ (no cross-links), integrate $\Psi^2 = a^2/(1-ab')^2$ from $b'=0$ to~$b$:
\[
\GFcomb\big|_{d=0}
= a + \int_0^b \frac{a^2}{(1-ab')^2}\,db'
= a + \frac{a^2 b}{1-ab}
= \frac{a}{1-ab},
\]
recovering the linear-polymer species GF.
\emph{Step~2}: at given~$b$, integrate the cross-link contribution from $d'=0$ to~$d$:
\begin{equation}\label{eq:xl_species_integral}
\GFcomb(a,b,d) = \frac{a}{1-ab} + \frac{a^2}{2}\int_0^d u(d')^4\,dd',
\end{equation}
where $u(d')$ is the physical root of the cubic~\eqref{eq:xl_cubic} at $(a,b,d')$.
The integral cannot be expressed in elementary functions for generic~$d$, but the formal power series in~$d$ is determined to all orders by the cubic.

\paragraph{Conservation equation.}
The total dimensionless atom count is
\begin{equation}\label{eq:xl_conservation}
N = a\,\PDeriv{\GFcomb}{a},
\end{equation}
relating the total concentration $c^{\mathrm{total}} = \gamma N$ to the free-monomer fugacity~$a$ at fixed bond variables $b,d$.

\paragraph{Singularity analysis and gelation.}
The coupled system~\eqref{eq:xl_reduced_system} has the form $(\Psi,\Psi_r) = \mathbf{F}(a,\Psi,\Psi_r)$ with
\[
F_1 = a(1+b\Psi)(1+d\Psi_r),
\qquad
F_2 = a(1+b\Psi)^2.
\]
By the Drmota--Lalley--Woods criterion~\eqref{eq:dlw_condition}, the dominant singularity occurs where $\det(I - \partial\mathbf{F}/\partial(\Psi,\Psi_r))=0$.
Computing the Jacobian:
\[
I - \begin{pmatrix}
ab(1+d\Psi_r) & ad(1+b\Psi)\\
2ab(1+b\Psi) & 0
\end{pmatrix},
\]
whose determinant, after substituting $\Psi_r = au^2$ and $1+d\Psi_r = 1+dau^2$, is
\begin{equation}\label{eq:xl_det}
\det(I-J) = 1 - ab - 3a^2 b d\,u^2.
\end{equation}
Setting $\det(I-J)=0$ gives $u^2 = (1-ab)/(3a^2 bd)$.
Substituting into the cubic~\eqref{eq:xl_cubic} and simplifying yields
$u_c = 3/(2(1-ab))$, and the \emph{gelation surface}:
\begin{equation}\label{eq:xl_gelation}
\boxed{27\,a^2 b d = 4(1-ab)^3.}
\end{equation}
This defines the critical fugacity $a_c(b,d)$ implicitly.
At $d=0$, the condition degenerates to $(1-ab)^3=0$, i.e.\ $a_c=1/b$ (the simple pole of the linear-polymer system).
For any $d>0$, the gelation point $a_c(b,d)<1/b$ is strictly less than the linear-polymer singularity: cross-links pull the gel point inward.
Because the system is nonlinear and irreducible, the Drmota--Lalley--Woods theorem guarantees a square-root branch point with
\[
[a^n]\,\GFcomb \sim C\,a_c^{-n}\,n^{-3/2},
\]
in contrast to the simple-pole behavior $[a^n]\,\GFcomb = b^{n-1}$ of linear polymers at $d=0$.
The critical concentration is $c_{\max} = \gamma\,a_c\,\GFcomb'(a_c)$, the maximum total atom concentration for which a pre-gel equilibrium exists.

\subsection{Cyclic case}
\label{sec:crosslinked_cyclic}

\paragraph{Linearized PDEs.}
Applying the cyclic linearization~\eqref{eq:linearized_cycle_identity} to each bond type, with $\widetilde{\GFensemble} = \exp(\GFspecies)$:
\begin{equation}\label{eq:xl_linearized_b}
\PDeriv{\widetilde{\GFensemble}}{b} = \frac{\partial^2 \widetilde{\GFensemble}}{\partial p\,\partial q}
\qquad\text{(heterotypic bond $b$)},
\end{equation}
\begin{equation}\label{eq:xl_linearized_d}
\PDeriv{\widetilde{\GFensemble}}{d} = \frac{1}{2}\,\frac{\partial^2 \widetilde{\GFensemble}}{\partial r^2}
\qquad\text{(homotypic bond $d$)}.
\end{equation}
Both are closed linear PDEs in~$\widetilde{\GFensemble}$.

\paragraph{Operator exponential.}
Define the operators $D_b = \partial^2/(\partial p\,\partial q)$ and $D_d = \partial^2/\partial r^2$.
Since $D_b$ acts on $(p,q)$ while $D_d$ acts on~$r$, the operators commute: $[D_b,D_d]=0$.
The joint solution~\eqref{eq:aggregate_operator_exponential} is therefore
\begin{equation}\label{eq:xl_operator_exp}
\widetilde{\GFensemble}
= \exp\!\left(b\,D_b + \tfrac{d}{2}\,D_d\right)\widetilde{\GFensemble}_0
= \exp(b\,D_b)\,\exp\!\left(\tfrac{d}{2}\,D_d\right)\widetilde{\GFensemble}_0,
\end{equation}
where the initial data~\eqref{eq:Z0_explicit} is $\widetilde{\GFensemble}_0 = \exp(a\,p\,q\,r)$.

\paragraph{Canonical partition function.}
By Theorem~\ref{thm:canonical_cyclic}, with $N$ atoms and $M_p=M_q=M_r=N$:
\begin{equation}\label{eq:xl_canonical}
\GFensemble_N
= \left.\exp\!\left(b\,D_b + \tfrac{d}{2}\,D_d\right)
p^N q^N r^N\right|_{p=q=r=1}.
\end{equation}
Because the operators act on disjoint site variables, the partition function \emph{factorizes}:
\begin{equation}\label{eq:xl_factorization}
\GFensemble_N = \GFensemble_N^{(b)}\;\cdot\;\GFensemble_N^{(d)},
\end{equation}
where
\begin{equation}\label{eq:xl_ZNb}
\GFensemble_N^{(b)}
= \sum_{i=0}^{N}\frac{b^{\,i}}{i!}\,(N)_i^2
= \sum_{i=0}^{N}\frac{N!}{(N-i)!}\binom{N}{i}\,b^{\,i}
\end{equation}
is the linear-polymer-with-rings partition function~\eqref{eq:rings_explicit}, and
\begin{equation}\label{eq:xl_ZNd}
\GFensemble_N^{(d)}
= \sum_{j=0}^{\lfloor N/2\rfloor}\frac{d^{\,j}}{2^j\,j!}\,(N)_{2j}
= \sum_{j=0}^{\lfloor N/2\rfloor}\binom{N}{2j}\,(2j-1)!!\;d^{\,j}
\end{equation}
is the \emph{matching polynomial}, where $(2j{-}1)!! = 1\cdot 3\cdot 5\cdots(2j{-}1)$ and $(-1)!!=1$.
The matching polynomial counts the number of ways to choose $j$ unordered $r$--$r$ pairs from $N$ monomers, weighted by~$d^j$.
The factorization reflects the independence of site pools: $b$-bonds consume $(p,q)$-sites while $d$-bonds consume $r$-sites, so the two bond types draw from disjoint resources.

\paragraph{Ensemble and species generating functions.}
The grand-canonical ensemble GF is
\begin{equation}\label{eq:xl_ensemble}
\GFensemble(a,b,d)
= \sum_{N=0}^{\infty}\frac{a^N}{N!}\,\GFensemble_N^{(b)}\,\GFensemble_N^{(d)},
\end{equation}
and the species GF is $\GFcomb = \log\GFensemble$.

\paragraph{Singularity analysis and gelation.}
The radius of convergence $R$ of $\GFcomb(a)$ (with $b,d$ as parameters) is determined by the growth rate of $\GFensemble_N^{(b)}\,\GFensemble_N^{(d)}/N!$.
The saddle-point equation~\eqref{eq:saddle_point_equation} for the cyclic system reads
\begin{equation}\label{eq:xl_saddle}
a\,\GFcomb'(a) = N,
\end{equation}
and the gelation threshold is $a_c = R$, the value at which~\eqref{eq:xl_saddle} ceases to admit a solution.
Because the factorization~\eqref{eq:xl_factorization} couples the asymptotics of both factors, the critical fugacity $a_c$ in the cyclic regime is shifted relative to the acyclic gelation surface~\eqref{eq:xl_gelation}: cycles provide additional species that absorb monomers, raising the effective monomer capacity and moving the gel point outward.
However, the singularity remains a square-root branch point with $n^{-3/2}$ coefficient asymptotics---cycles shift the gelation point but preserve the universality class.

\section{Outlook}
\label{sec:outlook}

The framework developed here suggests several immediate extensions.

\begin{itemize}
  \item \textbf{Further branching systems.} The cross-linked polymer system of Section~\ref{sec:crosslinked} illustrates commuting operator exponentials when bond types act on disjoint site pools.
  The full version of this paper develops a trivalent bond system with three distinguishable sites ($a, b, c$), one atom type ($x$), and two bond types ($y$ binding $a$--$b$, $z$ binding $a$--$c$).
  This enables a clean species recursion $T = x(1+yT)(1+zT)$ for rooted trees, Lagrange inversion giving Catalan-number structure, and a systematic treatment of cycles via the linearized ensemble PDEs.
  The trivalent system also illustrates non-commuting operator exponentials when multiple bond types share sites, complementing the commuting case treated in Section~\ref{sec:crosslinked}.

  \item \textbf{Beyond solvable motifs.} The operator-exponential form \eqref{eq:operator_exponential_solution} isolates the combinatorial effect of adding a cyclic closure but does not by itself make computation easy. The saddle-point approximation of Section~\ref{sec:gf_saddle} provides one systematic route, with $O(1/N)$ error (Appendix~\ref{app:saddle}), and the singularity analysis of Section~\ref{sec:singularity} characterizes the phase boundary at which the saddle-point method breaks down. Extending this to multivariate saddle points for systems with several conserved atom types, computing systematic higher-order finite-size corrections from the cumulant expansion~\eqref{eq:saddle_point_error}, and developing numerical coefficient-extraction methods adapted to realistic contact maps remain natural next steps.

  \item \textbf{Richer constraints and validity.} Physical validity constraints (sterics, saturation, geometry) prune the class of admissible assemblies and can introduce long-range dependence. Incorporating such constraints while preserving computable generating-function identities remains open.

  \item \textbf{Foundational unification.} Several nearby formalisms (species, groupoid cardinality, master-equation generating functions) point to a unified categorical/combinatorial foundation for equilibrium assembly. Clarifying these relationships should improve both conceptual transparency and computational robustness.

  \item \textbf{Euler characteristic as admissibility criterion.} The Euler characteristic $\chi = 1 - \vartheta$ provides a minimal admissibility filter and generalizes naturally to higher-order bond systems. The notion of \emph{magnitude} of enriched categories \cite{leinster2008euler} is a far-reaching generalization of Euler characteristic whose investigation in the assembly context is an intriguing open question.
\end{itemize}

\appendix

\section{Closure-operator expansions}
\label{app:closure-expansions}

This appendix records explicit coefficient formulas for the backbone expansion used in the closure-operator recursion.

\subsection{Setup}
If $L=L(k;x,y)$ is the acyclic species GF built using a bond type $k$ between sites $x$ and $y$, then
\begin{equation}\label{eq:app_backbone_pde}
\PDeriv{L}{k} = \PDeriv{L}{x}\,\PDeriv{L}{y}.
\end{equation}
Write the $k$-series expansion $L = \sum_{n=0}^{\infty}\frac{k^n}{n!}\,L_n(x,y)$, $L_0 = A(x,y)$.
Substituting and equating coefficients:
\begin{equation}\label{eq:app_Ln_recursion}
L_{n+1}
=
\sum_{i=0}^{n}\binom{n}{i}\,
\frac{\partial L_i}{\partial x}\,
\frac{\partial L_{n-i}}{\partial y},
\qquad n\ge 0.
\end{equation}

\subsection{First coefficients}
\begin{align}
L_{1}
&=
\frac{\partial A}{\partial x}\,
\frac{\partial A}{\partial y},
\label{eq:app_L1}\\[0.5ex]
L_{2}
&=
\frac{\partial^2 A}{\partial x^2}\left(\frac{\partial A}{\partial y}\right)^2
+2\frac{\partial A}{\partial x}\frac{\partial A}{\partial y}\frac{\partial^2 A}{\partial x\partial y}
+\left(\frac{\partial A}{\partial x}\right)^2\frac{\partial^2 A}{\partial y^2},
\label{eq:app_L2}\\[0.5ex]
L_{3}
&=
\frac{\partial ^3A}{\partial x^3}\left(\frac{\partial A}{\partial y}\right)^3
+6\frac{\partial^2A}{\partial x^2}\left(\frac{\partial A}{\partial y}\right)^2\frac{\partial^2 A}{\partial x\partial y}
+6\frac{\partial A}{\partial x}\frac{\partial A}{\partial y}\left(\frac{\partial^2A}{\partial x\partial y}\right)^2\nonumber\\
&\quad
+3\frac{\partial A}{\partial x}\left(\frac{\partial A}{\partial y}\right)^2\frac{\partial^3A}{\partial x^2\partial y}
+6\frac{\partial A}{\partial x}\frac{\partial^2 A}{\partial x^2}\frac{\partial A}{\partial y}\frac{\partial^2A}{\partial y^2}\nonumber\\
&\quad
+6\left(\frac{\partial A}{\partial x}\right)^2\frac{\partial^2A}{\partial x\partial y}\frac{\partial^2A}{\partial y^2}
+3\left(\frac{\partial A}{\partial x}\right)^2\frac{\partial A}{\partial y}\frac{\partial^3A}{\partial x\partial y^2}
+\left(\frac{\partial A}{\partial x}\right)^3\frac{\partial^3A}{\partial y^3}.
\label{eq:app_L3}
\end{align}

For linear polymers, $A(x,y)=Axy$, so all higher mixed derivatives vanish and the expansion recovers the closed form $L=Axy/(1-kA)$.

\section{Canonical factorial moments: proofs}
\label{app:proofs}

\begin{proof}[Proof of Theorem~\ref{thm:canonical_factorial_moments}]
Introduce a formal variable $\xi_s$ for each species $s$ and define the canonical probability generating function
\[
G_{\mathbf{N}}(\boldsymbol{\xi})
=\frac{\prod_a N_a!}{\GFensemble_{\mathbf{N}}}\;
[\mathbf{z}^{\mathbf{N}}]\;
\exp\!\Bigl(\sum_{s}\xi_s\,w_s\,\mathbf{z}^{\boldsymbol{\alpha}(s)}\Bigr).
\]
At $\boldsymbol{\xi}=\mathbf{1}$, $G_{\mathbf{N}}(\mathbf{1})=1$.
The $k$-th factorial moment of $n_s$ equals the $k$-th derivative of $G_{\mathbf{N}}$ at $\boldsymbol{\xi}=\mathbf{1}$.
Apply the shift $\boldsymbol{\xi}\to\boldsymbol{\xi}+\mathbf{1}$:
\[
G_{\mathbf{N}}(\boldsymbol{\xi}+\mathbf{1})
=\frac{\prod_a N_a!}{\GFensemble_{\mathbf{N}}}\;
[\mathbf{z}^{\mathbf{N}}]\;\GFensemble(\mathbf{z})\;\exp\!\Bigl(\sum_{s}\xi_s\,w_s\,\mathbf{z}^{\boldsymbol{\alpha}(s)}\Bigr).
\]
Taking $k_j$ derivatives in $\xi_{s_j}$ pulls down $w_{s_j}^{k_j}\,\mathbf{z}^{k_j\boldsymbol{\alpha}(s_j)}$.
Setting $\boldsymbol{\xi}=\mathbf{0}$:
\[
\Bigl\langle \prod_j(n_{s_j})_{k_j}\Bigr\rangle_{\!\mathbf{N}}
=\frac{\prod_a N_a!}{\GFensemble_{\mathbf{N}}}\;\prod_j w_{s_j}^{k_j}\;
[\mathbf{z}^{\mathbf{N}-\mathbf{m}}]\;\GFensemble(\mathbf{z})
=\prod_j w_{s_j}^{k_j}\;\prod_a(N_a)_{m_a}\;
\frac{\GFensemble_{\mathbf{N}-\mathbf{m}}}{\GFensemble_{\mathbf{N}}},
\]
using $[\mathbf{z}^{\mathbf{N}-\mathbf{m}}]\,\GFensemble = \GFensemble_{\mathbf{N}-\mathbf{m}}/\prod_a(N_a-m_a)!$ and simplifying the factorial ratio to the falling factorial $(N_a)_{m_a}$.
\end{proof}

For compact statements, let $\mathbf{k}=(k_s)_{s\in\Species}\in\Nat^{(\Species)}$ be a finitely supported multi-index and define
\[
(\mathbf{n})_{\mathbf{k}} := \prod_{s\in\Species}(n_s)_{k_s},\qquad
\mathbf{w}^{\mathbf{k}} := \prod_{s\in\Species}w_s^{\,k_s},\qquad
\boldsymbol{\alpha}(\mathbf{k}) := \sum_{s\in\Species}k_s\,\boldsymbol{\alpha}(s).
\]

\begin{proposition}[Multi-index form]
\label{prop:canonical_factorial_moments_multiindex}
For every finitely supported $\mathbf{k}\in\Nat^{(\Species)}$,
\begin{equation}\label{eq:canonical_factorial_moments_multiindex}
\langle (\mathbf{n})_{\mathbf{k}}\rangle_{\mathbf{N}}
=
\mathbf{w}^{\mathbf{k}}
\;\cdot\;
(\mathbf{N})_{\boldsymbol{\alpha}(\mathbf{k})}
\;\cdot\;
\frac{\GFensemble_{\mathbf{N}-\boldsymbol{\alpha}(\mathbf{k})}}{\GFensemble_{\mathbf{N}}}.
\end{equation}
\end{proposition}

\begin{proof}
Group equal species in Theorem~\ref{thm:canonical_factorial_moments}, writing $k_s$ for the total number of factors at species~$s$.
\end{proof}

\begin{proposition}[Telescoping depletion factorization]
\label{prop:canonical_factorial_moment_telescoping}
For $\mathbf{k}=\mathbf{k}^{(1)}+\mathbf{k}^{(2)}$:
\begin{equation}\label{eq:canonical_factorial_telescoping}
\langle(\mathbf{n})_{\mathbf{k}}\rangle_{\mathbf{N}}
=
\langle(\mathbf{n})_{\mathbf{k}^{(1)}}\rangle_{\mathbf{N}}
\cdot
\langle(\mathbf{n})_{\mathbf{k}^{(2)}}\rangle_{\mathbf{N}-\boldsymbol{\alpha}(\mathbf{k}^{(1)})}.
\end{equation}
In particular, for a fixed species $s$ and $k\ge 1$:
\begin{equation}\label{eq:canonical_factorial_single_species_product}
\langle (n_s)_k\rangle_{\mathbf{N}}
=
\prod_{j=0}^{k-1}\langle n_s\rangle_{\mathbf{N}-j\,\boldsymbol{\alpha}(s)}.
\end{equation}
\end{proposition}

\begin{proof}
Apply \eqref{eq:canonical_factorial_moments_multiindex} to $\mathbf{k}$, $\mathbf{k}^{(1)}$, and $\mathbf{k}^{(2)}$, using
$\mathbf{w}^{\mathbf{k}}
=
\mathbf{w}^{\mathbf{k}^{(1)}}
\mathbf{w}^{\mathbf{k}^{(2)}}$,
$(\mathbf{N})_{\boldsymbol{\alpha}(\mathbf{k})}
=
(\mathbf{N})_{\boldsymbol{\alpha}(\mathbf{k}^{(1)})}
(\mathbf{N}-\boldsymbol{\alpha}(\mathbf{k}^{(1)}))_{\boldsymbol{\alpha}(\mathbf{k}^{(2)})}$,
and the telescoping ratio
$\GFensemble_{\mathbf{N}-\boldsymbol{\alpha}(\mathbf{k})}/\GFensemble_{\mathbf{N}}
=
(\GFensemble_{\mathbf{N}-\boldsymbol{\alpha}(\mathbf{k}^{(1)})}/\GFensemble_{\mathbf{N}})
\cdot
(\GFensemble_{\mathbf{N}-\boldsymbol{\alpha}(\mathbf{k})}/\GFensemble_{\mathbf{N}-\boldsymbol{\alpha}(\mathbf{k}^{(1)})})$.
\end{proof}

\begin{proposition}[Recursive identities from partition-function ratios]
\label{prop:canonical_factorial_moment_recursions}
For any finitely supported $\mathbf{k}$:
\begin{equation}\label{eq:canonical_ratio_from_moment}
\frac{\GFensemble_{\mathbf{N}-\boldsymbol{\alpha}(\mathbf{k})}}{\GFensemble_{\mathbf{N}}}
=
\frac{\langle (\mathbf{n})_{\mathbf{k}}\rangle_{\mathbf{N}}}
{\mathbf{w}^{\mathbf{k}}\,(\mathbf{N})_{\boldsymbol{\alpha}(\mathbf{k})}}.
\end{equation}
For any species $s$:
\begin{equation}\label{eq:canonical_partition_recursion_from_mean}
\GFensemble_{\mathbf{N}}
=
\frac{w_s\,(\mathbf{N})_{\boldsymbol{\alpha}(s)}}{\langle n_s\rangle_{\mathbf{N}}}\;
\GFensemble_{\mathbf{N}-\boldsymbol{\alpha}(s)}.
\end{equation}
\end{proposition}

\begin{proof}
Equation \eqref{eq:canonical_ratio_from_moment} is \eqref{eq:canonical_factorial_moments_multiindex} solved for the ratio; \eqref{eq:canonical_partition_recursion_from_mean} follows by setting $\mathbf{k}=\mathbf{e}_s$.
\end{proof}

Thus higher-order factorial moments are recursively decomposed into lower-order ones at successively depleted compositions.
In particular, \eqref{eq:canonical_factorial_single_species_product} expresses $\langle (n_s)_k\rangle_{\mathbf{N}}$ as a product of first-order depleted-state factors.

\section{Saddle-point derivation and admissibility}
\label{app:saddle}

This appendix derives the saddle-point approximation stated in Proposition~\ref{prop:saddle_point_general}, computes the leading error correction~\eqref{eq:saddle_point_error}, and spells out the Flajolet--Sedgewick conditions that guarantee its validity.

\subsection{Derivation of the saddle-point formula}

We apply Theorem~VIII.4 of Flajolet and Sedgewick~\cite{flajolet2009analytic} (Hayman's admissibility theorem) to $f(z)=e^{\GFcomb(z)}$.
The auxiliary functions are $a(r) = r\,\GFcomb'(r)$ and $b(r) = r^2\,\GFcomb''(r) + r\,\GFcomb'(r)$.
The saddle point $x^*$ satisfies $a(x^*)=N$, i.e.\ $x^*\GFcomb'(x^*)=N$, and Hayman's theorem gives
\[
[z^N]\,e^{\GFcomb(z)}
\;\sim\;
\frac{e^{\GFcomb(x^*)}}{(x^*)^N\,\sqrt{2\pi\,b(x^*)}}.
\]
To connect with the Cauchy integral, parametrize the circle of radius $x^*$ by $z = x^*e^{i\theta}$ and define $h(z) = \GFcomb(z) - N\log z$, so that $h'(x^*)=0$.
The variance parameter at the saddle point is
\[
\lambda_N := -\left.\frac{d^2}{d\theta^2}h(x^*e^{i\theta})\right|_{\theta=0}
= (x^*)^2\,\GFcomb''(x^*) + N = b(x^*),
\]
which coincides with~\eqref{eq:lambda_def}.
The Gaussian approximation of the contour integral, justified by the H-admissibility conditions below, gives
\begin{equation}\label{eq:app_gaussian_eval}
[z^N]\,e^{\GFcomb(z)}
\;\sim\;
\frac{e^{\GFcomb(x^*)}}{(x^*)^N\,\sqrt{2\pi\lambda_N}}.
\end{equation}
Multiplying by $N!$ and applying Stirling's approximation $N!\sim \sqrt{2\pi N}\,(N/e)^N$:
\[
\GFensemble_N
= N!\,[z^N]\,e^{\GFcomb(z)}
\;\sim\;
\sqrt{\frac{N}{\lambda_N}}\;\exp\!\bigl(\GFcomb(x^*)\bigr)\,
\left(\frac{N}{e\,x^*}\right)^{\!N},
\]
which is~\eqref{eq:saddle_point_general}.

\subsection{Higher-order correction}

Write $u=i\theta$ and define the coefficients
\[
\kappa_j
:= \left.\frac{d^j}{du^j}h(x^*e^{u})\right|_{u=0},
\qquad j\ge 2,
\]
so that $\kappa_1=0$ (saddle-point equation), $\kappa_2 = (x^*)^2\GFcomb''(x^*)+N = \lambda_N$, and
\begin{align*}
\kappa_3 &= (x^*)^3\GFcomb'''(x^*) + 3(x^*)^2\GFcomb''(x^*) + N,\\
\kappa_4 &= (x^*)^4\GFcomb^{(4)}(x^*) + 6(x^*)^3\GFcomb'''(x^*) + 7(x^*)^2\GFcomb''(x^*) + N.
\end{align*}
(In general $\kappa_j = \sum_{k=1}^{j}S(j,k)\,(x^*)^k\GFcomb^{(k)}(x^*)$, where $S(j,k)$ are Stirling numbers of the second kind; the $k=1$ term gives $S(j,1)\,x^*\GFcomb'(x^*) = N$.)
The Taylor expansion of $h$ in $\theta$ is
\[
h(x^*e^{i\theta})
= h(x^*) - \tfrac{1}{2}\kappa_2\,\theta^2
- \tfrac{i}{6}\kappa_3\,\theta^3
+ \tfrac{1}{24}\kappa_4\,\theta^4
+ \cdots\,.
\]
Write the integrand as $e^{h(x^*)}e^{-\kappa_2\theta^2/2}\bigl(1 + \delta(\theta)\bigr)$ where $\delta$ collects the cubic and quartic corrections:
\[
\delta(\theta) = -\frac{i\kappa_3}{6}\theta^3 + \frac{\kappa_4}{24}\theta^4 - \frac{\kappa_3^2}{72}\theta^6 + \cdots\,.
\]
The odd power $\theta^3$ contributes only at second order (via $\delta^2$); the even powers contribute at first order.
Evaluating the Gaussian moments $\langle \theta^{2k}\rangle_{\kappa_2} = (2k{-}1)!!/\kappa_2^k$ gives:
\[
\frac{\int e^{h(x^*e^{i\theta})}d\theta}
     {\int e^{h(x^*)-\kappa_2\theta^2/2}d\theta}
= 1
+ \frac{1}{8}\,\frac{\kappa_4}{\kappa_2^2}
- \frac{5}{24}\,\frac{\kappa_3^2}{\kappa_2^3}
+ O\!\left(\kappa_2^{-2}\right).
\]
The $\kappa_4$ term comes from $\langle\theta^4\rangle = 3/\kappa_2^2$; the $\kappa_3^2$ term from $\langle\theta^6\rangle = 15/\kappa_2^3$.
Including the next Stirling correction $N! = \sqrt{2\pi N}(N/e)^N(1+1/(12N)+\cdots)$:
\begin{equation}\label{eq:app_error_expansion}
\frac{\GFensemble_N}{\widehat{\GFensemble}_N}
= 1
+ \frac{1}{8}\,\frac{\kappa_4}{\kappa_2^2}
- \frac{5}{24}\,\frac{\kappa_3^2}{\kappa_2^3}
+ \frac{1}{12N}
+ O\!\left(N^{-2}\right).
\end{equation}
Since $\kappa_j=\Theta(N)$ for all $j\ge 2$, each correction term is $O(1/N)$.

\subsection{H-admissibility: Flajolet--Sedgewick conditions}
\label{app:saddle_admissibility}

The rigorous justification for extending the Gaussian integral to $(-\infty,\infty)$ in~\eqref{eq:app_gaussian_eval}---and hence for the saddle-point formula---is provided by the theory of H-admissible (Hayman-admissible) functions~\cite{hayman1956generalisation}, as codified by Flajolet and Sedgewick~\cite{flajolet2009analytic} (Chapter~VIII).

\begin{definition}[H-admissibility {\cite[Definition~VIII.1]{flajolet2009analytic}}]
\label{def:h_admissible}
A function $f(z)=\sum_{n\ge 0}f_n\,z^n$ with non-negative coefficients and radius of convergence $R\in(0,\infty]$ is \emph{H-admissible} if $f(r)>0$ for $r$ in a neighborhood of~$R$ and, writing
\[
a(r) := r\,\frac{f'(r)}{f(r)},\qquad
b(r) := r\,a'(r) = r^2\,\frac{f''(r)}{f(r)} + a(r) - a(r)^2,
\]
the following three conditions hold:
\begin{itemize}
\item[\textbf{H1.}] $a(r)\to+\infty$ and $b(r)\to+\infty$ as $r\to R^-$.
\item[\textbf{H2.}] (Locality.) There exists $\delta(r)>0$ with $0<\delta(r)<\pi$ such that, uniformly for $|\theta|\le\delta(r)$ as $r\to R^-$,
\[
f(r\,e^{i\theta})\;\sim\; f(r)\,\exp\!\bigl(i\theta\,a(r) - \tfrac{1}{2}\theta^2\,b(r)\bigr).
\]
\item[\textbf{H3.}] (Capture.) Uniformly for $\delta(r)\le|\theta|\le\pi$ as $r\to R^-$,
\[
f(r\,e^{i\theta}) = o\!\left(\frac{f(r)}{\sqrt{b(r)}}\right).
\]
\end{itemize}
\end{definition}

\noindent
Condition \textbf{H2} ensures that the saddle-point neighborhood dominates the integral, while \textbf{H3} ensures the tail is negligible.
Together they justify the Gaussian approximation.

\begin{theorem}[Hayman's theorem {\cite[Theorem~VIII.4]{flajolet2009analytic}}]
\label{thm:hayman}
If $f(z)$ is H-admissible with $a(r_n)=n$ (defining the saddle point~$r_n$), then
\[
[z^n]\,f(z) \;\sim\; \frac{f(r_n)}{r_n^{\,n}\,\sqrt{2\pi\,b(r_n)}}
\qquad\text{as } n\to\infty.
\]
\end{theorem}

\subsection{Verification for assembly partition functions}
\label{app:saddle_verification}

\begin{proposition}
\label{prop:ensemble_admissible}
Let $\GFcomb(z)=\sum_{s\in\Species}g(s)$ be a species generating function with non-negative coefficients (which holds since each $g(s) = (1/\varphi(s))\prod z_b^{\beta(s,b)}\ge 0$) and suppose $\GFcomb$ has finite radius of convergence $R$ with $\GFcomb(r)\to\infty$ as $r\to R^-$.
Then $\GFensemble(z)=\exp(\GFcomb(z))$ is H-admissible.
\end{proposition}

\begin{proof}
Since the coefficients of $\GFcomb$ are non-negative, so are those of $\GFensemble = \exp(\GFcomb)$.
We verify the three conditions for $f = \GFensemble$:

The auxiliary functions are $a(r) = r\,\GFcomb'(r)$ and $b(r) = r^2\,\GFcomb''(r) + r\,\GFcomb'(r)$.

\textbf{H1.} Because $\GFcomb(r)\to\infty$ as $r\to R^-$ and $\GFcomb$ has non-negative coefficients, $r\,\GFcomb'(r)\to\infty$.
The second derivative $\GFcomb''(r)\ge 0$ (all coefficients non-negative), so $b(r)\ge a(r)\to\infty$.

\textbf{H2--H3.} These follow from the general closure theorem: if $g(z)$ is any power series with non-negative coefficients, $g(0)=0$, and $b_g(r)\to\infty$ as $r\to R^-$, then $\exp(g(z))$ is H-admissible~\cite[Theorem~VIII.5]{flajolet2009analytic}.
In our setting $g=\GFcomb$ satisfies these hypotheses.
\end{proof}

Applying Theorem~\ref{thm:hayman} to $\GFensemble$ and combining with Stirling's formula for the $N!$ prefactor in~\eqref{eq:cauchy_integral} recovers the saddle-point formula of Proposition~\ref{prop:saddle_point_general}.

\subsection{Range of validity}

The saddle-point equation $x^*\GFcomb'(x^*)=N$ has a solution $x^*\in(0,R)$ if and only if $N<\lim_{r\to R^-}r\,\GFcomb'(r)$.
For generating functions with finite radius $R$ (such as $\GFcomb = a/(1-ab)$ with $R=1/b$), this imposes an upper bound $N_{\max}(b)$ on the atom number for which the approximation applies.
Physically, $N_{\max}$ corresponds to the saturation or gelation threshold: beyond it, the mass-balance equation has no solution, signaling a phase transition.

For generating functions with infinite radius of convergence (e.g., entire functions arising from assembly systems with finitely many species), $\lim_{r\to\infty}r\,\GFcomb'(r)=\infty$ and the saddle-point approximation is valid for all $N$.

In either case, the relative error is $O(1/N)$ by~\eqref{eq:app_error_expansion}, with explicit dependence on the cumulants $\kappa_3,\kappa_4$ at the saddle point.
For the linear polymer examples of Sections~\ref{sec:cyclefree_linear} and~\ref{sec:cycles_linear}, the error is sub-percent by $N\approx 50$ (see the worked examples in the respective sections).

\section{Foundations: species, cardinality, stuff types}
\label{app:foundations}

This appendix provides the formal backbone for the symmetry weights $1/\varphi(s)$, the exponential/multiset passage from species to ensembles, and the appearance of differential operators in stochastic equilibrium.

\subsection{Combinatorial species}

A \emph{combinatorial species} is a functor $F:\mathbf{FinBij}\to \mathbf{Set}$, where $\mathbf{FinBij}$ is the groupoid of finite sets and bijections.
For each finite set $U$, the set $F[U]$ is the set of $F$-structures on $U$.
The associated exponential generating function is
\begin{equation}\label{eq:app_species_egf}
F(x)=\sum_{n\ge 0}\frac{|F[n]|}{n!}\,x^n.
\end{equation}
The division by $n!$ expresses counting labeled realizations modulo relabeling.
When structures have internal symmetries, the coefficient becomes a symmetry-weighted quotient by automorphisms---the origin of the $1/\varphi(s)$ weights in \eqref{eq:species_monomial}.
See \cite{joyal1981theorie,bergeron1998combinatorial} for the standard development.

\subsection{Groupoid cardinality}

\begin{definition}
Let $\mathcal{G}$ be a groupoid with finite automorphism groups.
The \emph{groupoid cardinality} is
$|\mathcal{G}| := \sum_{[x]\in \pi_0(\mathcal{G})}1/|\mathrm{Aut}(x)|$.
\end{definition}

In the paper's notation, $\varphi(s) = |\mathrm{Aut}(s)|$.
The species generating function $\GFspecies=\sum_{s} g(s)$ is a symmetry-weighted sum over isomorphism classes, with each term weighted by $1/\varphi(s)$ as in the groupoid-cardinality framework of \cite{baez2001finite}.

\subsection{Stuff types and the exponential formula}

A \emph{stuff type} is a functor $p:\mathcal{X}\to \mathbf{FinBij}$, where $\mathcal{X}$ is a groupoid.
The symmetry-corrected count of $U$-structures is the groupoid cardinality of the fiber $\mathcal{X}_U$.
The class $\mathrm{SET}(F)$ of finite sets (bags) of $F$-components has generating function
\begin{equation}\label{eq:app_exponential_formula}
\mathrm{SET}(F)(x) = \exp(F(x)),
\end{equation}
because choosing, for each $k\ge 1$, how many $k$-labeled components appear (independently) and dividing by the permutations of identical components yields factorial denominators.
In the paper, $\GFcomb$ is the component series and $\GFensemble = \exp(\GFcomb)$ is the multiset/ensemble series.

\subsection{Stochastic mechanics dictionary}

A state is $x:S\to \mathbb{N}$ with finite support.
The probability generating function is $G_t(z) = \sum_x p_t(x)\prod_{s} z_s^{x(s)}$.
Creating one copy of $s$ corresponds to multiplying by $z_s$; annihilating one corresponds to $\partial/\partial z_s$.
The master equation becomes a linear PDE in $G_t$ \cite{baez2018quantum,baez2015quantum}.
At equilibrium under detailed balance, the product-form measure has factorial denominators $x(s)!$ matching the combinatorics of multisets---this is why $\GFensemble = \exp(\GFcomb)$ appears both as a combinatorial ``bags of components'' construction and as an equilibrium partition function.

\section*{Acknowledgements}

I gratefully acknowledge the support of the Santa Fe Institute, where this work was largely developed during an Omidyar Postdoctoral Fellowship.
I also acknowledge the support of the Fontana Lab at Harvard Medical School.

\bibliographystyle{plainnat}
\bibliography{refs}

\end{document}